\documentclass[12pt]{article}
\usepackage{amssymb}
 \usepackage{latexsym}
 \usepackage{theorem}

 \theoremstyle{plain} \theoremstyle{}
{\theorembodyfont{\itshape}
\newtheorem{thm}{Theorem}[section]}
\theoremstyle{plain} \theoremstyle{} {\theorembodyfont{\itshape}
\newtheorem{lem}[thm]{Lemma}}
\theoremstyle{plain} \theoremstyle{} {\theorembodyfont{\itshape}
\newtheorem{prop}[thm]{Proposition}}
\theoremstyle{plain} \theoremstyle{} {\theorembodyfont{\itshape}
}
\theoremstyle{plain} \theoremstyle{} {\theorembodyfont{\itshape}
\newtheorem{cor}[thm]{Corollary}}
\theoremstyle{plain} \theoremstyle{} {\theorembodyfont{\itshape}
}
\theoremstyle{plain} \theoremstyle{} {\theorembodyfont{\itshape}
}
\theoremstyle{plain} \theoremstyle{} {\theorembodyfont{\rmfamily}
}
\theoremstyle{plain} \theoremstyle{} {\theorembodyfont{\rmfamily}
}
\theoremstyle{plain} \theoremstyle{} {\theorembodyfont{\rmfamily}
}
\theoremstyle{plain} \theoremstyle{} {\theorembodyfont{\rmfamily}
}
\theoremstyle{plain} \theoremstyle{} {\theorembodyfont{\rmfamily}
}
\theoremstyle{plain} \theoremstyle{} {\theorembodyfont{\rmfamily}
\newtheorem{eg}[thm]{Example}}
\theoremstyle{plain} \theoremstyle{} {\theorembodyfont{\rmfamily}
}
\theoremstyle{plain} \theoremstyle{} {\theorembodyfont{\rmfamily}
\newtheorem{egs}[thm]{Examples}}
\theoremstyle{plain} \theoremstyle{} {\theorembodyfont{\rmfamily}
\newtheorem{df}[thm]{Definition}}

\theoremstyle{plain} \theoremstyle{} {\theorembodyfont{\itshape}
}
\theoremstyle{plain} \theoremstyle{} {\theorembodyfont{\itshape}
}
\theoremstyle{plain} \theoremstyle{} {\theorembodyfont{\itshape}
}
\theoremstyle{plain} \theoremstyle{} {\theorembodyfont{\rmfamily}
}

\begin{document}
\newcommand{\lam}{\lambda}
\newcommand{\da}{\downarrow}
\newcommand{\Da}{\Downarrow}
\newcommand{\D}{\Delta}
\newcommand{\ua}{\uparrow}
\newcommand{\ra}{\rightarrow}
\newcommand{\la}{\leftarrow}
\newcommand{\Lra}{\Longrightarrow}
\newcommand{\Lla}{\Longleftarrow}
\newcommand{\rat}{\rightarrowtail}
\newcommand{\up}{\upsilon}
\newcommand{\Up}{\Upsilon}
\newcommand{\ep}{\epsilon}
\newcommand{\ga}{\gamma}
\newcommand{\Ga}{\Gamma}
\newcommand{\Lam}{\Lambda}
\newcommand{\CF}{{\cal F}}
\newcommand{\CG}{{\cal G}}
\newcommand{\CH}{{\cal H}}
\newcommand{\CI}{{\cal I}}
\newcommand{\CB}{{\cal B}}
\newcommand{\CT}{{\cal T}}
\newcommand{\CS}{{\cal S}}
\newcommand{\CV}{{\cal V}}
\newcommand{\CP}{{\cal P}}
\newcommand{\CQ}{{\cal Q}}
\newcommand{\cu}{{\underline{\cup}}}
\newcommand{\ca}{{\underline{\cap}}}
\newcommand{\nb}{{\rm int}}
\newcommand{\Si}{\Sigma}
\newcommand{\si}{\sigma}
\newcommand{\Om}{\Omega}
\newcommand{\bm}{\bibitem}
\newcommand{\bv}{\bigvee}
\newcommand{\bw}{\bigwedge}
\newcommand{\lra}{\longrightarrow}
\newcommand{\tl}{\triangleleft}
\newcommand{\tr}{\triangleright}
\newcommand{\y}{{\bf y}}
\newcommand{\colim}{{\rm colim}}
\title{{\Large Many-Valued Complete
Distributivity}\thanks{This work is supported by NCET and 973
Programs (No. 2002cb312200) of China.}}
\author{ {\normalsize Hongliang Lai,  Dexue Zhang}\\
{\small Yangtze Center of Mathematics, Sichuan University, Chengdu
610064, China} \\
{\small  E-mail: honglianglai@163.com, dxzhang@scu.edu.cn}}
\date{}
\maketitle\begin{abstract}Suppose $(\Om, *, I)$ is a commutative,
unital quantale. Categories enriched over $\Om$ can be studied as
generalized, or many-valued, ordered structures. Because many
concepts, such as complete distributivity, in lattice theory can
be characterized by existence of certain adjunctions, they can be
reformulated in the many-valued setting in terms of categorical
postulations. So, it is possible, by aid of categorical
machineries, to establish theories of many-valued complete
lattices, many-valued completely distributive lattices, and so on.
This paper presents a systematical investigation of many-valued
complete distributivity, including the topics: (1) subalgebras and
quotient algebras of many-valued completely distributive lattices;
(2) categories of (left adjoint) functors; and (3) the
relationship between many-valued complete distributivity and
properties of the quantale $\Om$. The results show that enriched
category theory is a very useful tool in the study of many-valued
versions of order-related mathematical entities.


\end{abstract}
\parskip 3pt

\section{Introduction}

{\bf Categories enriched over a quantale as many-valued ordered
structures.} Partially ordered sets are important structures in
logic, mathematics, and theoretical computer sciences
\cite{AJ,AbVi,Birkhoff, G80, Rosenthal90}. From the viewpoint of
category theory, a partially ordered set, and generally a
preordered set, is a special kind of categories, i.e., a category
with hom-sets either empty or a singleton. Following Lawvere
\cite{Law73, Law86}, this fact can be put in a different way. Let
$\{0, 1\}$ denote the complete lattice consisting of two elements
with the ordering $0\leq 1$. Then, ${\bf 2} = (\{0, 1\}, \wedge,
1)$ is a symmetric, monoidal, closed category. Enriched categories
over $\bf 2$ are just the preordered sets. So, preordered sets can
be investigated by aid of the categorical machinery. A nice
example is the wide use of Galois connections, a special case of
adjoint functors, in the theory of partially ordered sets. More
importantly, many interesting concepts in the theory of partially
ordered sets can be postulated as certain categorical properties!
For example,

(1) a partially ordered set $P$ is a lattice if and only if the
diagonal $P\lra P\times P$ has both a left adjoint and a right
adjoint;

(2) a partially ordered set $P$ is a complete lattice if and only
if the Yoneda embedding $\y:P\lra {\cal D}(P)$ has  a left
adjoint, where, ${\cal D}(P)$ denotes the partially ordered set of
lower sets in $P$ with the inclusion ordering and $\y(p) =\
\downarrow\!p= \{q\in P\ |\ q\leq p\}$; and

(3) a complete lattice $P$ is (constructive) completely
distributive if and only if the left adjoint of the Yoneda
embedding also has a left adjoint.

The categorical postulations of lattices, complete lattices, and
completely distributive lattices, can be easily reformulated for
categories enriched over any arbitrary symmetric, monoidal, closed
category. The aim of this paper is to study complete
distributivity for categories enriched over a complete, symmetric,
and monoidal closed small category, i.e., categories enriched over
a commutative, unital quantale $(\Om, *, I)$.

The study of quantale-enriched categories as generalized ordered
structures originated from the theory of {\it quantitative domain
theory}, see, e.g. \cite{BBR,FSW,Ru96,Sch03,Smyth87,Wa94,Wa97}.
The idea is as follows. Suppose that $A$ is a category enriched
over a commutative, unital quantale $\Om$. Then, for any two
elements $x, y\in A$, the enrichment $A(x, y)$, an element in
$\Om$, can be interpreted as the degree that $x$ precedes $y$, or,
the degree that $x$ is smaller than or equal to $y$. Therefore, an
$\Om$-category can be regarded as a {\it quantitative preordered
set}, in which the relation between two points is expressed by an
element in the quantale $\Om$ (regarded as the set of
truth-values), as opposed to the traditional {\it qualitative},
yes-or-no, relation in a preordered set.

{\bf Logic aspect of categories enriched over a quantale.} As just
mentioned in the above paragraph, in the study of categories
enriched over a commutative, unital quantale $\Om$, the quantale
can be regarded as the set of truth values. So, the study of
quantale-enriched categories has a strong logical flavor. This
aspect of enriched categories was emphasized by Lawvere early in
1973 in \cite{Law73} as {\it generalized pure logic}. The idea is
roughly as follows. Since $\Om$ is a monoidal closed category, for
each $\alpha\in \Om$, the functor $\alpha*(_-):\Om\lra\Om$ has a
right adjoint $\alpha\ra(_-):\Om\lra\Om$. That means, for all
$\alpha, \beta, \gamma\in\Om$, $$\alpha*\beta\leq\gamma\iff
\beta\leq \alpha\ra\gamma.$$ Therefore, if we interpret $\alpha,
\beta,$ and $\gamma$ as truth values, then the operation $*$ plays
a similar role as the logic connective {\it conjunction} and $\ra$
as the connective {\it implication}. The least element $0\in\Om$
can be regarded as the logical value {\it absurdity} and the unit
element $I$ as the value {\it true}. Moreover, for every set $X$,
a function $\lam: X\lra\Om$ can be regarded as a predicate on $X$,
the element $\lam(x)$ is the degree that $x$ has certain
attribute. And it is natural to interpret $\bw_{x\in X}\lam(x)$
and $\bv_{x\in X}\lam(x)$ as the truth degree for the logical
formulas $\forall x\lam$ and $\exists x\lam$ respectively. These
observations relate the study of quantale-enriched categories to
many-valued logic.

By a "many-valued logic" we mean a logic of which the truth-value
set is just a commutative, unital quantale. Such a logic is called
a {\it Monoidal Logic} in \cite{Ho95,EG01}. When the quantale
$\Om$ is a $BL$-algebra, this kind of logic has been extensively
investigated under the name {\it Basic Logic} in the literature,
see, e.g. H\'{a}jek \cite{Ha98,Ha99}; and when $\Om$ is a
commutative Girard quantale, such a logic is a commutative version
of {\it Linear Logic} initiated by Girard, see, e.g.
\cite{Girard,Rosenthal90,Yetter}. And, if $\Om$ is at the same
time a $BL$-algebra and a commutative Girard quantale, then $\Om$
must be an $MV$-algebra \cite{Ho95},  in this case we come back to
the {\it Many-Valued Logic} initiated by {\L}ukasiewicz
\cite{CDM}. However, by abuse of language, we shall call a
"monoidal logic" simply a "many-valued logic" in this paper
because that the truth value set is not a boolean algebra in
general and hence contains more than two elements. And so the
title of this paper.

{\bf Relationship to constructive complete distributivity.}  It is
well-known that there already exists a notion of lattice in any
topos \cite{MLM,W}. Since many concepts in lattice theory can be
postulated as categorical properties, they can be easily
reformulated in any topos. So, we can establish theories of
complete lattices, completely distributive lattices, and so on,
within any topos. For complete distributivity, this has already
been done under the name {\it constructive complete
distributivity} in a series of papers \cite{FW,RW2,RW4,W}. These
theories are developed within the framework of the internal logic
in a topos.

But, many-valued complete distributivity could be regarded as, to
some extent, a mathematical theory developed within the framework
of an {\it observed} logic. This can be roughly explained as
follows. Suppose that we are working in {\bf Set}, even with the
Axiom of Choice allowed if you prefer. Contemplation over
principles of our reasoning about mathematical entities, or
generally the states of affairs around us, leads us to observe
that not only the  Boolean algebra $\bf 2$, but also any
$MV$-algebra, any $BL$-algebra, any commutative unital quantale,
posses sufficient structures to act as a {\it ruler} or a {\it
criterion} in our reasoning about "states of affairs in reality"
\footnote{This phrase is quoted from Dirk van Dalen: Logic and
Structure, 4th edition, Springer, 2004.}. Then, if we, tempted by
the fun to know, take  this kind of structures as our
criteria\footnote{Sometimes, this is necessary  and helpful as
exemplified in \cite{CDM,Ha98}.}, or truth values, we are led to
theories of many-valued logics. Many-valued complete
distributivity discussed in this paper is a mathematical theory
developed within the framework of an observed logic with a
commutative, unital quantale as the set of truth values. This
method could be generalized to establish many-valued versions of
other mathematical entities.

Since both the notion of many-valued complete distributivity and
that of constructive complete distributivity are postulated as a
certain categorical property, categorical methods play an
essential role in the study of both of these notions. And more
interestingly, many-valued complete distributivity and
constructive complete distributivity have similar properties. The
reader can compare the following Theorem 1.1 and Theorem 1.2.

\begin{thm} Suppose that $\Om$ is an integral commutative quantale.
Then the following conditions are equivalent:

$(1)$ $\Om$ is a Girard quantale, i.e., it satisfies the law of
double negation.

$(2)$ The dual of every  completely distributive $\Om$-lattice is
also a completely distributive $\Om$-lattice.

$(3)$ $\Om^{\rm op}$ is a many-valued Heyting algebra.\end{thm}

\begin{thm} {\rm \cite{RW2, W}} In any topos $\cal E$ the following
conditions are equivalent:

$(1)$ $\cal E$ is Boolean.

$(2)$ The dual of every  constructive completely distributive
lattice is constructive completely distributive.

$(3)$ The dual lattice of the truth value set is a Heyting
algebra.\end{thm}

However, the exact relationship between many-valued complete
distributivity and constructive complete distributivity still
awaits  further investigation.

{\bf Related works.} Categories enriched over a commutative,
unital quantale $\Om$ have received wide attention in the
literature since the publication of the pioneering paper of
Lawvere \cite{Law73}. The order aspect of $\Om$-categories leads
the theory of quantitative domains, see, e.g. \cite{BBR,FSW,Ru96,
Smyth87,Wa94,Wa97}. These works have a strong background in
theoretical computer sciences, so most of them are concerned with
certain kind of directed completeness of $\Om$-categories. A
systematical investigation of directed completeness for
$\Om$-categories is presented in \cite{LZ05}. And it should be
pointed out that the order aspect of quantale-enriched categories
has also been studied under the name {\it fuzzy order}, see, e.g.
\cite{Be02,Be04}.

Most of the results in Section 3 on characterization of
completeness of $\Om$-categories are special cases of the general
results for enriched categories in \cite{Bo,Ke}. Because of the
simplicity of $\Om$, these results appear in an extremely simple
form. The equivalence between cocomplete $\Om$-categories and
$\Om$-modules \cite{AbVi,JT84,Rosenthal96} was first observed by
Stubbe \cite{Stubbe04}.

The notion of completely distributive $\Om$-lattices was
introduced in \cite{Stubbe05} and \cite{Zh}. Stubbe has obtained
many interesting results about many-valued complete distributivity
in \cite{Stubbe05}. It is necessary to make clear the relationship
between \cite{Stubbe05} and this paper. At first, \cite{Stubbe05}
focuses on characterizations of many-valued complete
distributivity in the category of cocomplete $\Om$-categories and
cocontinuous $\Om$-functors. This paper is, with emphasis on the
order aspect of completely distributive $\Om$-lattices, dealing
with the category of many-valued completely distributive
$\Om$-lattices and complete $\Om$-lattice morphisms. Secondly, the
approach in \cite{Stubbe05} is comparatively more {\it
sophisticated}, it depends heavily on the computation techniques
on weighted limits and weighted colimits developed in
\cite{Stubbe05a} for quantaloids. But, the  approach in this paper
is quite {\it elementary}. Thirdly, the most important, except
propositions \ref{totally below relation is interpolative} and
\ref{lower sets is CD} (with different proofs), there is little
overlap between the results about many-valued completely
distributive lattices in the article \cite{Stubbe05} and this one.

And, it should be pointed out that in \cite{Stubbe05} $\Om$ is not
assumed to be commutative. One problem with the absence of
commutativity is that the dual of an $\Om$-category is not an
$\Om$-category in general (Example \ref{Omega should be
commutative}). Though many results and  proofs in this paper can
be improved to cope with the absence of commutativity of $\Om$ as
in \cite{Stubbe05}, we have assumed the commutativity of $\Om$ in
order to be succinct.

{\bf Summary of the contents.} This paper is devoted to a
systematical investigation of many valued complete distributivity.
The contents are arranged as follows.

In section 2, basic notions of $\Om$-categories and complete
$\Om$-lattices are recalled.

In section 3, some equivalent descriptions of complete
$\Om$-lattices are given. Of particular interest is the
observation by Stubbe \cite{Stubbe04} that complete $\Om$-lattices
are essentially $\Om$-modules. These characterizations of complete
$\Om$-lattices shall be often employed in the following sections.

Section 4 introduces the notion of many-valued completely
distributive lattices, i.e., completely distributive
$\Om$-lattices,  and discusses some basic properties of these
objects.

Section 5 focuses on the subalgebras and quotient algebras of a
completely distributive $\Om$-lattice. The main result says that
the subalgebras of a completely distributive $\Om$-lattice $A$
correspond bijectively to the cocontinuous closure operators on
$A$ and the quotient algebras of a completely distributive
$\Om$-lattice $A$ correspond bijectively to the cocontinuous
kernel operators on $A$.

Section 6 is an application of the results in Section 5. It is
proved that the category of left adjoints between completely
distributive $\Om$-lattices is completely distributive by showing
that it is a quotient algebra of some completely distributive
$\Om$-lattice.

The last section, Section 7, deals with the question that whether
the dual of a completely distributive $\Om$-lattice is also
completely distributive. The main result in this section is
Theorem 1.1 stated in the above. This result relates many-valued
complete distributivity closely to properties of the truth-value
set $\Om$.

\section{Complete $\Om$-lattices}

We refer to \cite{Bo,ML} for general category theory, to
\cite{Bo,Ke,Law73} for enriched category theory, and to
\cite{Birkhoff, G80} for lattice theory.

Let $\Om$ be a complete lattice. The greatest element of $\Om$ is
denoted $1$ and  the least element of $\Om$ is denoted $0$. For
$U\subset \Om$, write $\bv U$ for the least upper bound of $U$ and
$\bw U$ for the greatest lower bound of $U$. Particularly, $\bv
\emptyset=0$ and $\bw\emptyset=1$.

A commutative, unital quantale is a triple $(\Om, I, *)$,
abbreviated as $\Om$, where, $\Om$ is a complete lattice, $I$ is a
fixed element in $\Om$, and $\ast:\Om \times\Om\lra\Om$, called
the tensor, is a commutative, associative binary operation such
that (1) $*$ is monotone on each variable; (2) $I$ is a unit
element for $*$, i.e. $\alpha\ast I=\alpha$ for every
$\alpha\in\Om$; and (3) for each $\alpha\in\Om$, the monotone
function $\alpha\ast(_-):\Om\lra\Om$ has a right adjoint
$\alpha\ra(_-): \Om\lra \Om$. The resulting binary operation $\ra
:\Om\times\Om\lra\Om$ is called the residuation operator, or
implication, corresponding to the tensor $\ast$.

Throughout this paper, $(\Om, *, I)$ will always denote a
commutative, unital quantale. And when there will be no confusion
with respect to the tensor $*$ and the unit $I$,  we often write
simply $\Om$ instead of $(\Om, *, I)$. Some basic properties of
the tensor operator and residuation operator are collected in the
following, most of them can be found in many places, for instance,
\cite{Be02,Ha98,Ho95,Rosenthal90}.
\begin{prop}
$(1)\ 0\ast \alpha=0\ {\rm for\ all}\ \alpha\in\Om.$

$(2)\ \alpha\ra \beta=\bv \{\gamma:\alpha\ast \gamma\leq \beta\}.$

$(3)\ I\ra \alpha=\alpha;\ \  0\ra \alpha=1 \ {\rm for\ all}\
\alpha\in\Om.$

$(4)\ (\alpha\ra \beta)\ast(\beta\ra \gamma)\leq(\alpha\ra
\gamma). $

$(5)\ \alpha\ra(\beta\ra \gamma)=(\alpha\ast \beta)\ra \gamma =
\beta\ra(\alpha\ra \gamma). $

$(6)\ ((\alpha\ra \beta)\ra \beta)\ra \beta=\alpha\ra \beta. $

$(7)\ \alpha\ast\bv _{j\in J}\beta_j=\bv _{j\in J}\alpha\ast
\beta_j. $

$(8)\ (\bv_{j\in J}\alpha_j)\ra \beta=\bw _{j\in J}(\alpha_j\ra
\beta); \ \
 \alpha\ra(\bw _{j\in J}\beta_j)=\bw _{j\in J}(\alpha\ra \beta_j). $

$(9)\ \bw_{\gamma\in\Om}\Big((\gamma\ra \alpha)\ra(\gamma\ra
\beta)\Big)= \alpha\ra \beta.$

$(10)\ \bw_{\gamma\in\Om}\Big((\alpha\ra \gamma)\ra(\beta\ra
\gamma)\Big)=\beta\ra \alpha.$
\end{prop}

If the unit $I$ coincides with the greatest element $1$ in $\Om$,
$\Om$ is called an {\it integral commutative quantale}, or, a {\it
complete residuated lattice}. An integral commutative quantale
$\Om$ is called a {\it complete $BL$-algebra} \cite{Ha98} if it
satisfies

$(11)\ \alpha*(\alpha\ra \beta)=\alpha\wedge \beta$
(divisibility); and

$(12)\ (\alpha\ra \beta)\vee(\beta\ra \alpha)=1$ (prelinearity).

$BL$-algebras are the algebras for the {\it Basic Logic} developed
in \cite{Ha98}.

A commutative unital quantale $\Om$ is called a {\it commutative
Girard quantale} if it satisfies the {\it law of double negation}:

$(13)\ \alpha=(\alpha\ra 0)\ra 0$.

This definition of commutative Girard quantale is taken from
H\"{o}hle \cite{Ho95} and it is stronger than the definition in
\cite{Rosenthal90}. A commutative Girard quantale $\Om$ is
necessarily integral since
$$I = (I\ra 0)\ra 0= 0\ra 0 = 1.$$

Girard  quantales are closely related to the {\it Linear Logic}
developed by Girard \cite{Girard}.

A {\it complete $MV$-algebra} is a commutative unital quantale
which is simultaneously a $BL$-algebra and a  Girard quantale
\cite{Ho95}. For a nice exposition of $MV$-algebras and their role
in many-valued logic, we refer to the monograph \cite{CDM}.

A category enriched over a commutative quantale $\Om$, or an
$\Om$-category, is a pair $(A, {\rm hom})$ with $A$ a set and
${\rm hom}$ a function assigning to every ordered pair of
$(a,b)\in A\times A$ an element ${\rm hom}(a,b)\in \Om$, such that
(1) $I\leq {\rm hom}(a,a)$ for every $a\in A$ (reflexivity); and
(2) ${\rm hom}(a,b)\ast {\rm hom}(b,c)\leq {\rm hom}(a,c)$ for all
$a,b,c\in A$ (transitivity).

In an $\Om$-category  $(A, {\rm hom})$, $A$ is called the
underlying set of $(A, {\rm hom})$ and the function ${\rm hom}$ is
called the {\it hom functor}. We often write simply $A$ for an
$\Om$-category and $A(x,y)$ for hom$(x,y)$ if the hom functor is
clear from the context. And in this case we write $|A|$ for the
underlying set of $A$.

An $\Om$-functor between $\Om$-categories $A$ and $B$ is a
function $f:A\lra B$ such that $A(a,b)\leq B(f(a),f(b))$ for all
$a,b\in A$. An $\Om$-functor $f$ is called an {\it $\Om$-isometry}
if $A(a,b)= B(f(a),f(b))$ for all $a,b\in A$. If an $\Om$-isometry
$f$ is also bijective on the underlying sets, it will be called an
{\it $\Om$-isomorphism}, or an isomorphism for short.
$\Om$-functors are composed by composing the underlying functions
on sets.

An $\Om$-category $A$ can also be regarded as an $\Om$-valued
preordered set with the value $A(x, y)\in\Om$  being interpreted
as the degree to which $x$ is smaller than or equal to $y$. An
$\Om$-functor is also called  an {\it $\Om$-monotone function}
since the condition $A(a,b)\leq B(f(a),f(b))$ asserts that if
$a\leq b$ in $A$, then $f(a)\leq f(b)$ in $B$. In this paper we
switch freely between the terms of $\Om$-categories and
$\Om$-preordered sets, and between $\Om$-functors and
$\Om$-monotone functions. When we want to emphasize the
categorical aspect of $A$, we say $A$ is an $\Om$-category; and
when we want to emphasize the  order aspect of $A$, we say $A$ is
an $\Om$-preordered set. And so for the terms $\Om$-functors and
$\Om$-monotone functions.

Suppose $A$ is an $\Om$-category. We define a binary relation
$\leq$ on the underlying set of $A$ in the following way: $a\leq
b$ if $A(a,b)\geq I$. It is easily seen that $\leq$ is a preorder,
i.e. a reflexive and transitive relation, on $|A|$. For each
$\Om$-category, we write $A_0$ for the preordered set $(|A|,
\leq)$. In this way, we obtain a forgetful functor
$(_-)_0:\Om$-$\bf Cat$$\lra {\bf PrOrd}$ from the category
$\Om$-${\bf Cat}$ of $\Om$-categories to the category {\bf PrOrd}
of preordered sets.

Two elements $x$ and $y$ in an $\Om$-category $A$ are said to be
{\it isomorphic} if $A(x,y)\geq I$ and $A(y,x)\geq I$. An
$\Om$-category $A$ is called {\it antisymmetric} if different
elements in $A$ are always non-isomorphic, or equivalently, $A_0$
is a partially ordered set. An anti-symmetric $\Om$-category is
also called an {\it partially $\Om$-ordered set}.

In the following examples we list some  methods to construct
$\Om$-categories. These methods are somewhat standard in category
theory and it is hard to find where they appeared for the first
time, so, we don't include any reference here.

\begin{egs} (1) (The canonical $\Om$-category structure on $\Om$)
Let $\Om(\alpha, \beta) = \alpha\ra\beta$. Then, by Proposition
2.1, it is easy to check that $\Om$ is a partially $\Om$-ordered
set.

(2) (Discrete $\Om$-categories) Given a set $X$ and $x, y\in X$,
let $X(x, y) = I$ if $x=y$ and $X(x, y) =0$ if $x\not=y$. Then $X$
becomes an $\Om$-category. Such $\Om$-categories are called
discrete since that for any $\Om$-category $B$, every function
from $X$ to $B$ is an $\Om$-functor. We write $\bf 1$ for the
discrete $\Om$-category consisting of exactly one element.

(3) (Terminal object) Let $X =\{x\}$ be a singleton and $X(x, x) =
1$, the top element in $\Om$. Then $X$ is an $\Om$-category and it
is the terminal object in the category of $\Om$-categories, which
shall be denoted $\top$ in the sequel.

(4) (Dual $\Om$-category) Suppose $A$ is an $\Om$-category. Let
$A^{\rm op}(a,b)=A(b,a)$ for all $a,b\in A$. Then $A^{\rm op}$ is
also an $\Om$-category, called the dual of $A$.

(5) (Subcategory) Let $A$ be an $\Om$-category and $B$ is a subset
of $A$. For all $x, y\in B$, let $B(x, y) = A(x, y)$. Then $B$
becomes an $\Om$-category, called a {\it subcategory} of $A$.

(6) (Product category) Suppose $\{A_i:i\in J\}$ is a family of
$\Om$-categories, the product of $\{A_i:i\in J\}$ in the category
$\Om$-{\bf PrOrd} is given by
\[\prod_{i\in J}A_i(a,b)=\bw_{i\in J}A_{i\in J}(a_i,b_i),
\ a=(a_i)_{i \in I},\ b=(b_i)_{i\in J}.\]

(7) (Functor category) Given $\Om$-categories $A$ and $B$, denote
the set of all the $\Om$-functors from $A$ to $B$  by $[A,B]$. Let
$[A,B](f,g)=\bw_{x\in A}B(f(x),g(x))$ for all $f,g\in[A,B]$.  Then
$[A,B]$ becomes an $\Om$-category, called the {\it functor
category} from $A$ to $B$.

If $X$ is a discrete $\Om$-category, then $[X, B]$ consists of all
the function from $X$ to $B$. Particularly, given an
$\Om$-category $A$, let $|A|$ denote the $\Om$-category obtained
by equipping the underlying set of $A$ with the discrete
$\Om$-categorical structure. Sometimes, we write $[\Om^A]$
for$[|A|,\Om]$ in the sequel. Clearly, both $[A^{\rm op}, \Om]$
and $[A, \Om]$ are subcategories of $[\Om^A]$.
\end{egs}

\begin{eg}\label{Omega should be commutative} Suppose that
$(\Om, *, I)$ is a unital quantale, i.e., $(\Om, *, I)$ satisfies
the conditions for a commutative unital quantale except, possibly,
the commutativity of the binary operation $*$. For all $a,b\in
\Om$, let  $\Om(a,b)=a\ra b= \bv\{r\in\Om:a*r\leq b\}$. Then
$(\Om,\ra)$ becomes an $\Om$-category. However, the dual $\Om^{\rm
op}$ of $(\Om,\ra)$, given by $\Om^{\rm op}(a,b)=b\ra a$, is an
$\Om$-category if and only if $\Om$ is commutative. In fact, if
$\Om^{\rm op}$ is an $\Om$-category, then for all $a,b\in \Om$,
$$b*a\leq[a\ra (a*b)]*a=\Om^{\rm op}(a*b,a)*\Om^{\rm
op}(a,I)\leq\Om^{\rm op}(a*b,I)=a*b.$$ Exchanging the role of $a$
and $b$, we have that $a*b\leq b*a$. Thus, $\Om$ is
commutative.\end{eg}

\begin{df} Given an $\Om$-category $A$, the Yoneda embedding is
the function ${\bf y}:A\lra[A^{\rm op},\Om]$ given by ${\bf
y}(a)(x)=A(x,a)$ and the  co-Yoneda embedding is the function
${\bf y}':A\lra[A,\Om]^{\rm op}$ given by ${\bf y}'(a)(x)=A(a,x).$
\end{df}

The following lemma shows that both the Yoneda embedding and the
co-Yoneda embedding are $\Om$-isometries.

\begin{lem}\label{yoneda} {\rm (Yoneda)} {\rm(1)}  For all $a\in A$
and $\phi\in[A^{\rm op},\Om]$, $[A^{\rm op},\Om]({\bf
y}(a),\phi)=\phi(a)$.

{\rm(2)}  For all $a\in A$ and $\psi\in[A ,\Om]$, $[A ,\Om]({\bf
y}'(a),\psi)=\psi(a)$. \ \  $\Box$ \end{lem}

If $A$ is a (classical) preordered set. Then a $\bf 2$-functor
$\psi:A\lra{\bf 2}$ is, regarded as a characteristic function,
precisely an upper set of $A$ and a {\bf 2}-functor
$\phi\in[A^{\rm op},{\bf 2}]$ precisely a lower set of $A$. Thus,
an $\Om$-functor $\psi\in[A,\Om]$ shall also be called an upper
$\Om$-subset and $\phi\in[A^{\rm op},\Om]$ a lower $\Om$-subset in
$A$. Particularly, $\y(a) = A(-, a)$ and $\y'(a)= A(a, -)$ are
called respectively the principal lower $\Om$-subset and the
principal upper $\Om$-subset of $A$ generated by $a$.

\begin{prop}{\rm (\cite{LZ05})}  Let $A$ be an $\Om$-category.
Then for all $\alpha\in\Om$, $\CF\subset [A,\Om]$, and $\psi\in
[A,\Om]$, the followings hold:

$(1)$ $\bv\CF\in [A,\Om]$. Particularly, as the join of an empty
family, $0_A\in[A,\Om]$.

$(2)$ $\bw\CF\in[A,\Om]$. Particularly, as the meet of an empty
family, $1_A\in[A,\Om]$.

$(3)$ $\alpha\ast\psi\in [A,\Om]$.

$(4)$ $\alpha\ra\psi\in [A,\Om]$.

$(5)$ $\psi\ra \alpha\in[A^{\rm op}, \Om]$ and
$\psi=\bw_{\alpha\in \Om}\Big((\psi\ra \alpha)\ra \alpha\Big)$.\ \
$\Box$\end{prop}

For each $\mu\in\Om^A$ and $x\in A$, let $$(\ua\!\mu)(x)=\bv_{y\in
A}\mu(y)\ast A(y,x)=\bv_{y\in A}\mu(y)*{\bf y}'(y)(x),$$ and
$$(\da\!\mu)(x)=\bv_{y\in A}\mu(y)\ast A(x,y)=\bv_{y\in A}\mu(y)*{\bf
y}(y)(x).$$

It is easy to verify that $\ua\!\mu\in[A, \Om]$ and $\ua\!\mu$ is
the smallest upper $\Om$-subset in $A$ which is bigger than or
equal to $\mu$ under the pointwise order. Dually,  $\da\!\mu$ is
the smallest lower $\Om$-subset in $A$ which is bigger than or
equal to $\mu$ under the pointwise order. Particularly, $\phi\in
[A,\Om]\iff\ua\!\phi=\phi$ and $\phi\in [A^{\rm
op},\Om]\iff\da\!\phi=\phi.$ And we left it to the reader to check
that the two operations $\ua$ and $\da$ are $\Om$-functors, i.e.
$\ua: [\Om^A]\lra[A, \Om]$ and $\da: [\Om^A]\lra[A^{\rm op},\Om]$
are $\Om$-functors.

The following definition is a special case of the general concept
of enriched adjunctions in the theory of enriched categories
\cite{Bo,Ke}.

\begin{df}\cite{Wa94,Wa97} A pair of $\Om$-functors $f:A\lra B$ and
$g:B\lra A$ is said to be an $\Om$-adjunction (or simply an
adjunction) between $\Om$-categories $A$ and $B$, if
$B(f(a),b)=A(a,g(b))$ for all $a\in A$ and $b\in B$. In this case,
we say $f$ is a left adjoint of $g$ and $g$ is a right adjoint of
$f$. \end{df}

It is easily seen that if $(f, g)$ is an $\Om$-adjunction between
$\Om$-categories $A$ and $B$, then the pair of monotone functions
$f: A_0\lra B_0$ and $g: B_0\lra A_0$ is a Galois connection, i.e.
a $\bf 2$-adjunction, between the preordered set $A_0$ and $B_0$.

\begin{prop} \label{characterization of omega adjunction}
{\rm (\cite{LZ05})} Suppose $A$ and $B$ are $\Om$-categories and
$f:A\lra B,\ g:B\lra A$ are functions. Then  the following
conditions are equivalent:

{\rm(1)} $(f, g)$  is an $\Om$-adjunction.

{\rm(2)} $f$ is an $\Om$-functor and $B(f(a),b)=A(a,g(b))$ for all
$a\in A$ and $b\in B$.

{\rm(3)} $g$ is an $\Om$-functor and $B(f(a),b)=A(a,g(b))$ for all
$a\in A$ and $b\in B$.

These conditions imply

{\rm (4)} If $(f,g)$ is an $\Om$-adjunction between
$\Om$-categories $A$ and $B$, then $f\circ g\circ f=f$ and $g\circ
f\circ g=g$.\ \  $\Box$\end{prop}

\begin{prop}\label{injec is embed}
If $(f,g)$ is an $\Om$-adjunction between
$\Om$-categories $A$ and $B$, then

{\rm (1)} The following conditions are equivalent:

$(a)$ $f$ is injective,

$(b)$ $g\circ f={\rm id}_A$, and

$(c)$ $g$ is surjective.

In this case, $f$ is an isometry.

{\rm (2)} The following conditions are equivalent:

$(a')$ $f$ is surjective,

$(b')$ $f\circ g={\rm id}_B$, and

$(c')$ $g$ is injective.

In this case, $g$ is an isometry.\ \ $\Box$
\end{prop}

Let $f:A\lra B$ be an $\Om$-monotone function. For each
$\psi\in[B, \Om]$, let $f^\la(\psi)=\psi\circ f$. Then we obtain
an $\Om$-functor $f^\la:[B,\Om]\lra[A,\Om].$

\begin{prop}\label{kan extension} {\rm \cite{Law73}} The $\Om$-functor
$f^\la:[B,\Om]\lra[A,\Om]$ has both left and right adjoints. The
left adjoint is given by $f_\ell^\ra:[A, \Om]\lra[B, \Om]$, where
for all $\psi\in[A, \Om]$,
$$ f_\ell^\ra(\psi)(y)=\bv_{x\in A}\psi(x)*B(f(x),y).$$
The right adjoint is given by $f_r^\ra:[A, \Om]\lra[B, \Om]$,
where for all $\psi\in[A, \Om]$,$$ f_r^\ra(\psi)(y)=\bw_{x\in
A}B(y,f(x))\ra\psi(x).\ \ \Box$$
\end{prop}

Categorically, $f_\ell^\ra(\psi)$ and $f_r^\ra(\psi)$ are the
Kan-extensions of $\psi$ along $f$ \cite{Bo, Law73}.

Similarly, for each $\phi_B\in[B^{\rm op}, \Om]$ and
$\phi_A\in[A^{\rm op},\Om]$, let $f^\la(\phi_B)=\phi_B\circ
f\in[A^{\rm op}, \Om]$ and  $$f^\ra_\ell(\phi_A)(y)=\bv_{x\in
A}\phi_A(x)*B(y,f(x))$$ for all $y\in B$. Then both $f^\la:[B^{\rm
op},\Om]\lra[A^{\rm op},\Om]$ and $f^\ra_\ell:[A^{\rm
op},\Om]\lra[B^{\rm op},\Om]$ are $\Om$-functors.

Suppose that $f: A\lra B$ is an $\Om$-functor and  $\phi: A\lra
\Om$ is a function (not necessarily an $\Om$-functor). For each
$b\in B$, let
$$f(\phi)(b)=\bv_{x\in f^{-1}(b)}\phi(x).$$
Then $f(\phi)$ is called the image of $\phi$ under $f$. We left it
to the reader to check that for any $\psi\in[A, \Om]$,
$f^\ra_\ell(\psi) =\ \ua\!f(\psi)$ and for any $\phi\in[A^{\rm
op}, \Om]$, $f^\ra_\ell(\phi) =\ \da\!f(\phi)$.

\begin{eg}\label{lower is left adjoint}
Let $A$ be an $\Om$-category. If we write $i: [A, \Om]\lra[\Om^A]$
for the inclusion functor. Then $(\ua, i)$ is an $\Om$-adjunction.
Similarly, if we write $i: [A^{\rm op}, \Om]\lra[\Om^A]$ for the
inclusion functor. Then $(\da, i)$ is also an $\Om$-adjunction.

This fact can be put as a special case of the above proposition.
Let $|A|$ denote the underlying set of $A$ with the discrete
$\Om$-category structure. Then, the identity function ${\rm
id}:|A|\lra A$ is an $\Om$-functor. Clearly, ${\rm id}^\la:[A^{\rm
op},\Om]\lra[\Om^A]$ is exactly the inclusion $i: [A^{\rm op},
\Om]\lra[\Om^A]$ and ${\rm id}^\ra_\ell(\phi)\ =\ \da\!\phi$ for
all $\phi\in\Om^A$. Hence, $(\da, i)$ is an
$\Om$-adjunction.\end{eg}

\begin{df}\label{def of complete omega cats}{\rm \cite{LZ05,Wa94,Wa97}}
An $\Om$-category $A$ is called complete if the co-Yoneda
embedding $\y':A\lra[A,\Om]^{\rm op}$ has a right adjoint $\inf:
[A, \Om]^{\rm op}\lra A$. And $A$ is called cocomplete if the
Yoneda embedding $\y:A\lra[A^{\rm op},\Om]$ has a left adjoint
$\sup: [A^{\rm op}, \Om]\lra A$. \end{df}

Cocomplete $\Om$-categories are a special case of the {\it total}
(enriched) categories  in \cite{Ke86}. However, because of the
simplicity of $\Om$ and the fact that we don't have size problems
here, properties of these $\Om$-categories will become much
simpler and more elegant.

\begin{prop}{\rm \cite{Zh}} For an $\Om$-category $A$, the
following are equivalent:

$(1)$  $A$ is complete.

$(2)$ The composition functor $i\circ\y': A\lra[A, \Om]^{\rm
op}\lra[\Om^A]^{\rm op}$  has a right adjoint.

$(3)$  $A$ is cocomplete.

$(4)$ The composition functor $i\circ\y: A\lra[A^{\rm
op},\Om]\lra[\Om^A]$  has a left adjoint. \ \ \ $\Box$
\end{prop}

That $i\circ\y: A\lra[A^{\rm op},\Om]\lra[\Om^A]$ has a left
adjoint  amounts to say that for each $\phi\in[\Om^A]$, there is
an element $a\in A$ such that for all $x\in A$, $$A(a, x) =
[\Om^A](\phi, \y(x)) = \bw_{z\in A}(\phi(z)\ra A(z, x)).$$ The
condition that for all $x\in A$, $A(a, x) = \bw_{z\in
A}(\phi(z)\ra A(z, x))$, can be interpreted as the statement that
for all $x\in A$, $a$ is smaller than or equal to $x$ if and only
if $\phi$ is contained in the principal lower $\Om$-subset
generated by $x$. In other words, $a$ is the supremum of $\phi$ in
$A$ \cite{Wa94, Wa97}.

Generally, we say that a function $\phi:A\lra\Om$ has a supremum
if there is some element (unique up to isomorphisms) $a\in A$ such
that for all $x\in A$, $A(a, x) = \bw_{z\in A}(\phi(z)\ra A(z,
x))$. Then, by Proposition \ref{characterization of omega
adjunction}, an $\Om$-category $A$ is cocomplete if and only if
every function $\phi: A\lra\Om$ has a supremum in $A$.

If $A$ is cocomplete,  the left adjoint of $i\circ\y:A\lra[\Om^A]$
is given by $\sup\circ\da:[\Om^A]\lra[A^{\rm op}, \Om]\lra A$.
Thus, for each $\phi\in[\Om^A]$,  the supremum of $\phi$ is
$\sup(\da\!\phi)$. Thus, we shall write simply $\sup\phi$ for the
supremum of $\phi$ in the sequel.

Similarly, we can define the {\it infimum} $\inf\phi\
(=\inf(\ua\!\phi))$ of a function $\phi: A\lra\Om$ in $A$ and
obtain similar results.

\begin{eg}Both the singleton discrete $\Om$-category $\bf 1$ and
the terminal $\Om$-category $\top$ are complete.\end{eg}

\begin{prop} {\rm (\cite{Bo})} For each complete $\Om$-category
$A$, the underlying preordered set $A_0$ of $A$ is
complete.\end{prop}

A function  $\mu:A\lra \Om$ is said to be a finite $\Om$-subset if
the set $\{a\in A:\mu(a)\neq 0\}$ is finite.

\begin{df} An $\Om$-functor $f:A\lra B$ is said to preserve
(finite) sups if  $f(\sup{_A}\mu)=\sup{_B}f(\mu)$ whenever
$\sup{_A}\mu$ exists for any (finite) $\mu\in[\Om^A]$; and $f$ is
said to preserve (finite) infs if $f(\inf{_A}\mu)=\inf{_B}f(\mu)$
whenever $\inf{_A}\mu$ exists for any (finite) $\mu\in[\Om^A]$.
$f$ is also said to be continuous if it preserves infs and
cocontinuous if it preserves sups.
\end{df}

The following theorem is a generalization of the properties of
Galois connections between complete lattices. It relates the
existence of left (right) adjoints to the preservation of infima
(suprema resp.).

\begin{thm}\label{adjunction}{\rm(\cite{LZ05,Stubbe04})} Let
$f:A\lra B$ and $g:B\lra A$ are $\Om$-functors. If $(f,g)$ is an
$\Om$-adjunction, then $f: A\lra B$ preserves sups and $g: B\lra
A$ preserves infs. Conversely, if $A$ is cocomplete then $f$ has a
right adjoint whenever $f$ preserves sups; and  if $B$ is complete
then $g$ has a left adjoint whenever $g$ preserves infs. \ \ \
$\Box$
\end{thm}

\section{Equivalents of complete $\Om$-lattices}

\begin{df} {\rm (\cite{Bo,Ke})} An $\Om$-category $A$ is said
to be tensored if for all $\alpha\in\Om, x\in A$, there is an
element $\alpha\otimes x\in A$, called the tensor of $\alpha$ and
$x$, such that $A(\alpha\otimes x, y) = \alpha\ra A(x, y)$ for any
$y\in A$. $A$ is said to be cotensored if, for all $\alpha\in\Om$
and $y\in A$, there is some $\alpha\rightarrowtail x\in A$, called
the cotensor of $\alpha$ and $x$, such that $A(z,
\alpha\rightarrowtail x) = \alpha\ra A(z, x)$ for  any $z\in A$.
\end{df}

By definition, every complete (hence cocomplete) $\Om$-category is
both tensored and cotensored. And it is easily seen that the
tensor of $\alpha$ and $x$ in $A$ is just the cotensor of $\alpha$
and $x$ in $A^{\rm op}$. Since left adjoints preserve suprema,
they preserve tensors in the sense that $f(\alpha\otimes x) =
\alpha\otimes f(x)$, see Proposition \ref{characterization of left
adjoint} for a proof. Dually, right adjoints preserve cotensors.

\begin{egs} (1) $(\Om, \ra)$ is tensored and cotensored. Indeed,
$\alpha\otimes x = \alpha*x$ and $\alpha\rightarrowtail x =
\alpha\ra x$.

(2) For every $\Om$-category $A$, the functor category $[A^{\rm
op}, \Om]$ is tensored and cotensored. For all $\lam\in[A^{\rm
op}, \Om]$ and $\alpha\in\Om$,  the tensor $\alpha\otimes \lam$ of
$\alpha$ and $\lam$ in $[A^{\rm op}, \Om]$ is $\alpha*\lam$ and
the cotensor $\alpha\rightarrowtail \lam$ of $\alpha$ and $\lam$
in $[A^{\rm op}, \Om]$ is $\alpha\ra\lam$. Similarly, the tensor
and cotensor of $\alpha\in \Om$ and $\mu\in[A, \Om]$ in $[A, \Om]$
are given by $\alpha*\mu$ and $\alpha\ra\mu$ respectively.
\end{egs}

\begin{df}A complete $\Om$-lattice is an antisymmetric, complete
(hence cocomplete) $\Om$-category.\end{df}

\noindent{\bf Convention}: Suppose $A$ is a complete
$\Om$-lattice. Then $A_0$ is a complete lattice. Given a subset
$(x_t)_{t\in T}$ of $A$, the least upper bound of $(x_t)_{t\in T}$
in the complete lattice $A_0$ is called the {\it join} of
$(x_t)_{t\in T}$, $\bv_{t\in T}x_t$ in symbols; and the greatest
lower bound of $(x_t)_{t\in T}$ in  $A_0$ is called the {\it meet}
of $(x_t)_{t\in T}$, $\bw_{t\in T}x_t$ in symbols. We reserve the
notations $\sup$ and $\inf$ for supremum and infimum in $A$. That
is, for every function $\phi: A\lra\Om$, $\sup\phi$ stands for the
{\it supremum} of $\phi$ in $A$, and $\inf\phi$ for the {\it
infimum} of $\phi$ in $A$. So, for example, given a function
$\phi: A\lra\Om$, $\bv_{x\in A}\phi(x)$ denotes the least upper
bound of $\{\phi(x) \ |\ x\in A\}$ in the complete lattice $\Om$;
meanwhile $\sup\phi$ is an element in $A$, the supremum of $\phi$
in $A$.

\begin{prop}\label{tensor and cotensor} Suppose $A$ is a
tensored and cotensored $\Om$-category such that $A_0$ is a
complete lattice. Then, the tensor $\otimes: \Om\times A_0\lra
A_0$ and cotensor $\rightarrowtail:\Om\otimes A_0\lra A_0$ satisfy
the following conditions:

$(1)$ $A(\alpha\otimes x, y) = \alpha\ra A(x,y)=A(x,\alpha\rat
y)$. Hence $\alpha\otimes x\leq y \iff \alpha\leq A(x, y)\iff
x\leq\alpha\rat y$.

$(2)$ $({\rm i})\  I\otimes x=x;$ \ \  $({\rm i}')\ I\rat x=x $.

$(3)$ $({\rm ii})\ (\alpha*\beta)\otimes x =
\alpha\otimes(\beta\otimes x); $ \ \  $({\rm ii}')\
(\alpha*\beta)\rat x = \alpha\rat(\beta\rat x)$.

$(4)$ For any $x\in A$, the function $(_-)\otimes x:\Om\lra A_0$
is a left adjoint of the function $A(x,-): A_0\lra\Om$. Hence,
$$({\rm iii})\  \Big(\bv_{t\in T}\alpha_t\Big)\otimes x = \bv_{t\in
T}(\alpha_t\otimes x); \ \ ({\rm iii}')\ A\Big(x,\bw_{t\in
T}x_t\Big)=\bw_{t\in T}A(x,x_t).$$

$(5)$ For any $x\in A$, the function $(_-)\rat x:\Om\lra A_0^{\rm
op}$ is a left adjoint of the function $A(-,x): A_0^{\rm
op}\lra\Om$. Hence,
$$({\rm iv})\ A\Big(\bv_{t\in T}x_t,x\Big)=\bw_{t\in T}A(x_t,x);\
\ ({\rm iv}')\ \Big(\bv_{t\in T}\alpha_t\Big)\rat x=\bw_{t\in
T}(\alpha_t\rat x).$$

$(6)$ For any $\alpha\in\Om$, the function $\alpha\otimes(_-):
A_0\lra A_0$ is a left adjoint of the function $\alpha\rat(_-):
A_0\lra A_0$. Hence,
$$({\rm v})\ \alpha\otimes\Big(\bv_{t\in T}x_t\Big)=\bv_{t\in
T}\alpha\otimes x_t;\ \  ({\rm v}')\ \alpha\rat\Big(\bw_{t\in
T}x_t\Big)=\bw_{t\in T}(\alpha\rat x_t).$$
\end{prop}

\noindent{\bf Proof.} \ (1) and (2) follows from definition
immediately.

(3)  For any $x,y\in A$,
\begin{eqnarray*}A(\alpha\otimes(\beta\otimes x),y)
&=&\alpha\ra A(\beta\otimes x,y)=\alpha\ra(\beta\ra A(x,y))\\
&=&(\alpha*\beta)\ra A(x,y)= A((\alpha*\beta)\otimes x,y),
\end{eqnarray*} and \begin{eqnarray*}A(x,\alpha\rat(\beta\rat
y)&=&\alpha\ra A(x,\beta\rat x)=\alpha\ra(\beta\ra A(x,y))\\
&=&(\alpha*\beta)\ra A(x,y)=A(x,(\alpha*\beta)\rat
y).\end{eqnarray*}

(4), (5) and (6) follow from (1) straightforwardly. \ \ $\Box$

\begin{prop} Suppose $A$ and $B$ are complete $\Om$-lattices. Then
$f: A\lra B$ is an $\Om$-functor if and only if $(1)$ $f: A_0\lra
B_0$ preserves order and $(2)$ $\alpha\otimes f(x)\leq
f(\alpha\otimes x)$ for all $\alpha\in\Om$ and $x\in A$.\end{prop}

\noindent{\bf Proof.} \ This is because $f$ is an $\Om$-functor if
and only if for all $\alpha\in\Om, x, y\in A$, $\alpha\leq A(x,
y)$ implies that $\alpha\leq B(f(x), f(y))$. That means,
$\alpha\otimes x\leq y$ implies $\alpha\otimes f(x)\leq f(y)$,
which is equivalent to that $f: A_0\lra B_0$ preserves order and
$\alpha\otimes f(x)\leq f(\alpha\otimes x)$ for all $\alpha\in\Om$
and $x\in A$. \ \ \ $\Box$

\begin{prop} {\rm (\cite{Stubbe04})}
\label{characterization of left adjoint} Suppose $A$ and $B$ are
tensored $\Om$-categories, and $f: A\lra B$ is an $\Om$-functor.
Then, the followings are equivalent:

$(1)$ $f$ is a left adjoint.

$(2)$ $f: A_0\lra B_0$ is a left adjoint and $f$ preserves tensors
in the sense that $f(\alpha\otimes x)=\alpha\otimes f(x)$.
\end{prop}

\noindent{\bf Proof.} \ $(1) \Rightarrow(2)$: We need only show
that $f$ preserves tensors. Suppose $g:B\lra A$ is the right
adjoint of $f$. For all $\alpha\in\Om,x\in A,y\in B$,
\begin{eqnarray*}B(f(\alpha\otimes x),y)&=&A(\alpha\otimes x,g(y))
=\alpha\ra A(x,g(y))\\
&=&\alpha\ra B(f(x),y)=B(\alpha\otimes
x,y).\end{eqnarray*}Therefore, $f(\alpha\otimes x)=\alpha\otimes
f(x)$.

$(2) \Rightarrow(1)$: Suppose $g:B_0\lra A_0$ is a right adjoint
of $f:A_0\lra B_0$. For all $\alpha\in\Om,x\in A,y\in B$,
\begin{eqnarray*}\alpha\leq B(f(x),y)&\iff&\alpha\otimes f(x)\leq
y\iff f(\alpha\otimes x)\leq y\\
&\iff& \alpha\otimes x\leq g(y)\iff \alpha\leq
A(x,g(y)).\end{eqnarray*}

\ \ \ $\Box$

Similarly, we have the following.

\begin{prop} {\rm (\cite{Stubbe04})} Suppose $A$ and $B$ are
cotensored $\Om$-categories and $f: A\lra B$ is an $\Om$-functor.
Then, the followings are equivalent:

$(1)$ $f$ is a right adjoint.

$(2)$ $f: A_0\lra B_0$ is a right adjoint and $f$ preserves
cotensors in the sense that $f(\alpha\rightarrowtail
x)=\alpha\rightarrowtail f(x)$.\ \ \ $\Box$
\end{prop}

\begin{thm} {\rm (\cite{Stubbe04})} An antisymmetric, tensored and
cotensored $\Om$-category $A$ is a cocomplete (hence complete)
$\Om$-lattice if and only if $A_0$ is a complete lattice.
\end{thm}

\noindent{\bf Proof.} \ We need only check  the sufficiency.
Suppose that $A_0$ is complete. When regarded as a function
$A_0\lra [A^{\rm op}, \Om]_0$, the Yoneda embedding $\y$ is a
right adjoint by \ref{tensor and cotensor}(4);  and it preserves
cotensors by \ref{tensor and cotensor}(1). Thus, $\y:A\lra [A^{\rm
op}, \Om]$ has a left adjoint by the above proposition.\ \ $\Box$

The following proposition shows that the suprema and infima in a
complete $\Om$-lattice $A$ can be completely described by the
lattice structure of $A_0$ and the tensors and cotensors in $A$.

\begin{prop}Suppose $A$ is a complete $\Om$-lattice. We have the
following:

$(1)$ For every $\lam\in\Om^A$, $\sup\lam=\bv_{x\in
A}(\lam(x)\otimes x)$.

$(2)$ For every $\mu\in\Om^A$, $\inf\mu=\bw_{x\in
A}(\mu(x)\rightarrowtail x)$.\end{prop}

\noindent{\bf Proof.}\ \ We prove (2) for example. Because the
$\Om$-functor $\inf:[A, \Om]^{\rm op}\lra A$ is, by definition, a
right adjoint, it preserves meets and cotensors. And the meets and
cotensors in $[A,\Om]^{\rm op}$ are exactly the joins and tensors
in $[A,\Om]$ respectively. Therefore,
\begin{eqnarray*} \inf\mu=\inf(\ua\!\mu) & =&
\inf\Big(\bv_{x\in A}(\mu(x)*\y'(x))\Big)\\
& =& \bw_{x\in A}\inf(\mu(x)*\y'(x))\\
&=& \bw_{x\in A}(\mu(x)\rat\sup\circ\y(x))\\
&=&\bw_{x\in A}(\mu(x)\rat x),\end{eqnarray*} where the last
equality is from that $\sup\circ\y={\rm id}_A$ because $\y$ is
injective. \ \ $\Box$

The above results show that the structure of a complete
$\Om$-lattice $A$ can be completely described by the complete
lattice structure of $A_0$, the tensor $\otimes$ and the cotensor
$\rat$ on  $A$. By (i), (ii), (iii) and (v) in \ref{tensor and
cotensor}, a complete $\Om$-lattice $A$ is an $\Om$-module in the
category of complete lattices and join-preserving functions
\cite{AbVi,JT84,Rosenthal96}. Conversely, given an $\Om$-module in
category of complete lattices and join-preserving functions, i.e.,
a complete lattice $A_0$ and binary operation $\otimes: \Om\times
A_0\lra A_0$ which satisfies (i), (ii), (iii) and (v) in
\ref{tensor and cotensor}, let $A(x, y) = \bv\{\alpha\in\Om \ |\
\alpha\otimes x\leq y\}$. Then, (1) $A(x, x)\geq I$ for all $x\in
A$; (2) for all $x, y, z\in A$,
\begin{eqnarray*}A(x, y)*A(y, z)&=&\Big(\bv\{\alpha
\ |\ \alpha\otimes x\leq y\}\Big)*\Big(\bv\{\beta\ |\
\beta\otimes y\leq z\}\Big)\\
&=& \bv\{\alpha*\beta \ |\ \alpha\otimes x\leq y, \beta\otimes
y\leq z \}\\
&\leq& \bv\{\alpha \ |\ \alpha\otimes x\leq z \}\\
&=& A(x, z).\end{eqnarray*}

Thus, $A$ becomes an $\Om$-category. Moreover, we say that $A$ is
cocomplete. To this end, we show that for all $\lam: A^{\rm
op}\lra\Om$, the supremum of $\lam$ in $A$ is given by
$$\sup\lam = \bv_{x\in A}(\lam(x)\otimes x).$$

By definition of $A(x, y)$, for all $\alpha\in\Om$, $\alpha\leq
A(x, y)\iff\alpha\otimes x\leq y$.  Therefore,
\begin{eqnarray*} \alpha\leq A\Big({\bv}_{x\in A}(\lam(x)\otimes x),
y\Big) &\iff& \bv_{x\in
A}\alpha\otimes(\lam(x)\otimes x)\leq y\\
&\iff& \bv_{x\in
A}(\alpha*\lam(x))\otimes x\leq y\\
&\iff& \forall x\in A, \alpha*\lam(x) \leq A(x,y)\\
&\iff& \alpha\leq\bw_{x\in A}(\lam(x)\ra A(x,y)),
\end{eqnarray*} therefore, $$A\Big({\bv}_{x\in
A}(\lam(x)\otimes x), y\Big) = \bw_{x\in A}(\lam(x)\ra A(x,y)),$$
which means that $\sup\lam =\bv_{x\in A}(\lam(x)\otimes x)$.

Therefore, complete $\Om$-lattices and $\Om$-modules are
essentially the same things. This fact was first pointed out by
Stubbe in \cite{Stubbe04,Stubbe05}.

Dually, given a complete lattice $A_0$ and a binary operator
$\rat:\Om\otimes A_0^{\rm op}\lra A^{\rm op}$ which satisfies the
conditions ${\rm (i'),\ (ii'),\ (iv')\ and\ (v')}$ in \ref{tensor
and cotensor}, let $A(x, y) = \bv\{\alpha\in\Om \ |\ x\leq \alpha
\rat y\}$. Then $A$ becomes a complete $\Om$-lattice.

\begin{egs} (1) \cite{Bo,LZ05} Suppose that $A$
is an $\Om$-category. Then the $\Om$-category $[A^{\rm op},\Om]$
is tensored and cotensored. Since the underlying poset of $[A^{\rm
op},\Om]$ is a complete lattice, $[A^{\rm op},\Om]$ is a complete
$\Om$-lattice. For any function $G:[A^{\rm op},\Om]\lra\Om$,
$$\sup\ G=\bv_{\phi\in[A^{\rm op},\Om]}G(\phi)*\phi,\ \
\inf G=\bw_{\phi\in[A^{\rm op},\Om]}G(\phi)\ra\phi.$$

(2) (\cite{Bo}) The $\Om$-category $(\Om, \ra)$ is a complete
$\Om$-lattice since $\Om\cong[\Om^{\bf 1}]$, where ${\bf 1}$ is
the singleton discrete $\Om$-category. Therefore, for all
$\mu\in\Om^\Om$,
$$\inf\mu=\bw_{y\in\Om}\mu(y)\ra y,\ \ \ \
\sup\mu=\bv_{y\in\Om}\mu(y)\ast y.$$
\end{egs}

\begin{thm}{\rm (Tarski Fixed-point Theorem)} Suppose $A$ is a complete
$\Om$-lattice and $f:A\lra A$ is an $\Om$-functor. Then the set of
fixed points of $f$, ${\bf Fix}(f)=\{x\in A:f(x)=x\}$, as a
subcategory of $A$, is also a complete $\Om$-lattice.\end{thm}

\noindent{\bf Proof.} \ Firstly, we show that the set of prefixed
points of $f$, $M=\{x\in A:x\leq f(x)\}$, where $\leq$ is the
order on $A_0$, as a subcategory of $A$, is complete. For any
$B\subset M$, $\bv B\leq\bv_{b\in B}f(b)\leq f(\bv B)$, which
implies that $M$ is closed under the formation of arbitrary joins
in $A_0$, hence $M_0$ is a complete lattice. To see that $M$ is a
complete $\Om$-lattice, it is enough to show $M$ is closed under
the formation tensors. For all $\alpha\in\Om,\ x\in M$,
$\alpha\otimes x\leq\alpha\otimes f(x)\leq f(\alpha\otimes x)$ and
thus $\alpha\otimes x\in M$.

Secondly, note that the image of $M$ under $f$ is also contained
in $M$, thus we can restrict the domain and codomain of $f$ and
get a new $\Om$-functor $f':M\lra M$, which is also an
$\Om$-functor from $M^{\rm op}$ to $M^{\rm op}$. The prefixed
points of $f':M^{\rm op}\lra M^{\rm op}$ are exactly the fixed
points of $f$. Thus, the $\Om$-category ${\bf Fix}(f)^{\rm op}$,
as a subcategory of $M^{\rm op}$, is complete and then ${\bf
Fix}(f)$ is a complete $\Om$-lattice.\ \ \ $\Box$

\begin{prop}\label{functor category is complete} {\rm
(\cite{LZ05})} Suppose $B$ is complete $\Om$-category. Then for
any $\Om$-category $A$, the functor category $[A, B]$ is complete.
\end{prop}

\noindent{\bf Proof.} \ This conclusion has already been proved in
\cite{LZ05} by showing that all the weighted limits exist. Here we
include another, relatively simpler, proof here.

(1) $[A, B]_0$ is complete. At first, observe that $f\leq g$ in
$[A, B]_0$ if and only if $f(x)\leq g(x)$ in $B_0$ for all $x\in
A$. Now suppose $(f_t)_{t\in T}$ is a family of $\Om$-functors
from $A$ to $B$. Define $f(x) = \bv_{t\in T}f_t(x)$ for all $x\in
A$, where the join is taken in $B_0$. It suffices to show that
$f\in[A, B]$. Actually, for all $x, y\in A$, \begin{eqnarray*}
B(f(x), f(y)) &=& B\Big({\bv}_{t\in T}f_t(x), {\bv}_{t\in
T}f_t(y)\Big)\\
&=& \bw_{t\in T}B\Big(f_t(x), {\bv}_{t\in T}f_t(y)\Big)\\
&\geq&\bw_{t\in T}B(f_t(x), f_t(y))\\
&\geq& A(x, y).
\end{eqnarray*}

(2) Denote the tensor in $B$ by $\otimes_B$. For each
$\alpha\in\Om$, define a function $\alpha\otimes_B(_-): B\lra B$
by $\alpha\otimes_B(_-)(x) =\alpha\otimes_Bx$. Then
$\alpha\otimes_B(_-)$ is an $\Om$-functor since if $\beta\leq B(x,
y)$ then $\beta\otimes_Bx\leq y$, and hence $\beta\otimes_B(\alpha
\otimes_B x) = \alpha\otimes_B(\beta\otimes_B x) \leq
\alpha\otimes_B y$. Thus, $\beta\leq B(\alpha\otimes_Bx,
\alpha\otimes_By)$.

(3) Given $\alpha\in\Om, f\in[A, B]$, let $(\alpha\otimes f)(x) =
\alpha\otimes_B f(x)$. Then $\alpha\otimes f\in[A, B]$ since it is
the composition of $f$ and $\alpha\otimes_B(_-)$. We leave it to
the reader to check that $\otimes:\Om\times[A,B]_0\lra[A,B]_0$
satisfies the conditions (i), (ii), (iii) and (v) in \ref{tensor
and cotensor}, hence $[A,B]$ is a complete $\Om$-lattice. \ \
$\Box$

\begin{df} A $\Om$-functor $f:A\lra B$ between the
complete $\Om$-lattices is called a complete $\Om$-lattice
morphism if it has both left and right adjoints.\end{df}

Complete $\Om$-lattices and complete $\Om$-lattice morphisms form
a category, which shall be denoted $\Om$-{\bf CLat}.

\begin{prop}{\rm (\cite{LZ05})} Suppose $\{A_i:i\in J\}$ is a
family of complete $\Om$-lattices. The product $\prod_{i\in J}A_i$
is also complete.\end{prop}

\noindent{\bf Proof.}  \  Clearly, $\prod_{i\in J}(A_i)_0$ is a
complete lattice. We define an action of $\Om$ on $\prod_{i\in
J}(A_i)_0$ by $\alpha\otimes (a_i)_{i\in J}=(\alpha\otimes_i
a_i)_{i\in J}$, $\alpha\in\Om$, $(a_i)_{i\in J}\in \prod_{i\in
J}(A_i)_0$, where $\otimes_i$ is the tensor on  $A_i$. It is easy
to check that $\otimes$ satisfies the conditions (i), (ii), (iii)
and (v) in \ref{tensor and cotensor}, so it determines an
$\Om$-categorical structure on $\prod_{i\in J}(A_i)_0$. What
remains is to show that $\prod_{i\in J}(A_i)_0$ together with this
$\Om$-categorical structure coincides with the product
$\prod_{i\in J}A_i$. In fact, for any $\alpha\in \Om$ and $a,b\in
\prod_{i\in J}A_i$,
\begin{eqnarray*}\alpha\otimes a\leq b&\iff&\forall i\in J,
\alpha\otimes_i a_i\leq b_i\\
&\iff& \forall i\in J, \alpha\leq A_i(a_i,b_i)\\
&\iff &\alpha\leq\bw_{i\in J}A_i(a_i,b_i)=\prod_{i\in
J}A_i(a,b).\end{eqnarray*}\ \  $\Box$

For each $j\in J$, the projection $p_j:\prod_{i\in J}A_i\lra A_j$
is a complete $\Om$-lattice morphism. The left adjoint is given by
$$q_j:A_j\cong A_j\times\prod_{i\neq j}\{0_{A_i}\}
\hookrightarrow\prod_{i\in J}A_i,$$ and the right adjoint is given
by
$$q_j':A_j\cong A_j\times\prod_{i\neq j}\{1_{A_i}\}\hookrightarrow
\prod_{i\in J}A_i.$$ Thus, $\prod_{i\in J}A_i$ is the product of
$\{A_i:i\in J\}$ in the category $\Om$-{\bf CLat}.

\begin{prop} \label{equalizer exist} Suppose $f$ and $g$ are
complete $\Om$-lattice morphisms from a complete $\Om$-lattice $A$
to a complete $\Om$-lattice $B$. Then the equalizer of $f$ and $g$
exists in the category $\Om$-{\bf CLat}.\end{prop}

\noindent{\bf Proof.}\ \   Let $E=\{x\in A:f(x)=g(x)\}$ and
$E(x,y)=A(x,y)$. Then $E$ is a subcategory of $A$. It is enough to
show that $E$ is complete and the embedding $i:E\lra A$ is a
complete $\Om$-lattice morphism. Take any $\Om$-subset $\mu:E\lra
\Om$. We have that $f(\sup_Ai(\mu))=\sup_Bf\circ
i(\mu)=\sup_Bg\circ i(\mu)=g(\sup_Ai(\mu))$ because $f$ and $g$
preserves sups. Thus $\sup_Ai(\mu)\in E$ and then
$\sup_E\mu=\sup_Ai(\mu)$. Similarly, we can check that
$\inf_E\mu=\inf_Ai(\mu)$. Therefor, $E$ is complete and $i$ is a
complete $\Om$-lattice morphism indeed.\ \ $\Box$

Recall that a category is complete if and only if the products and
the equalizers exist. Clearly, $\Om$-{\bf CLat} is a complete
category.

\section{Completely distributive $\Om$-lattices}

\begin{df}\cite{Stubbe05,Zh} A complete $\Om$-lattice $A$ is said to be
completely distributive if the functor $\sup: [A^{\rm op},\Om]\lra
A$ has a left adjoint, denoted, $\Downarrow: A\lra [A^{\rm
op},\Om]$.
\end{df}

Clearly, when $\Om=\bf 2$, completely distributive $\Om$-lattices
coincide with the constructive completely distributive lattices in
\cite{FW,RW2,RW4,W}

\begin{eg} The $\Om$-category $(\Om, \ra)$ is completely distributive,
i.e., the $\Om$-functor $\sup:[\Om^{\rm op},\Om]\lra\Om$ has a
left adjoint. This is a special case of the general result
Proposition \ref{lower sets is CD} below. However, we shall
construct here a left adjoint of $\sup:[\Om^{\rm op},\Om]\lra\Om$
explicitly. At first, for each $\phi\in[\Om^{\rm op}, \Om]$, we
have:

(1) $\phi:\Om\lra\Om$ is a decreasing function;

(2) For all $x\in\Om, x*\phi(x)\leq \phi(I)$ since $x=I\ra
x\leq\phi(x)\ra \phi(I)$;

(3) For all $x\in \Om$, $\phi(I)\leq(x\ra I)\ra\phi(x)$ since
$x\ra I\leq \phi(I)\ra\phi(x)$.

Thus, for each $\phi\in[\Om^{\rm op}, \Om]$, $$\sup\phi =
\bv_{x\in\Om}x*\phi(x)=\phi(I)\ \ \ {\rm and}\ \ \
\bw_{x\in\Om}((x\ra I)\ra\phi(x))=\phi(I).$$

For each $x\in\Om$, let $\Da(x):\Om\lra\Om$ be given by $\Da(x)(t)
= x*(t\ra I)$ for all $t\in\Om$. Clearly $\Da(x)\in[\Om^{\rm op},
\Om]$ and we claim that $\Da:\Om\lra[\Om^{\rm op}, \Om]$ is a left
adjoint of $\sup:[\Om^{\rm op}, \Om]\lra \Om$. In fact, for all
$\phi\in[\Om^{\rm op}, \Om]$ and $x\in\Om$,
\begin{eqnarray*} [\Om^{\rm op}, \Om](\Da(x), \phi)
&=&\bw_{t\in\Om}((x*(t\ra I))\ra\phi(t)) = \bw_{t\in\Om}x\ra((t\ra
I)\ra\phi(t))\\
&=&x\ra\bw_{t\in\Om}(t\ra I)\ra\phi(t) = x\ra\phi(I)\\
& =&\Om(x, \sup\phi).\end{eqnarray*}\end{eg}

Suppose $A$ is a completely distributive $\Om$-lattice. It is
easily seen that for all $a\in A$ and $\lam\in[A^{\rm op}, \Om]$,
$\sup(\Downarrow\!(a)) =a$  and $\sup\lam\geq a$ if and only if
$\Downarrow\!(a)\leq \lam$. And, following the terminologies in
the series \cite{FW,RW2,RW4,W} on constructive completely
distributive lattices, we call $\Downarrow\!(a)(x)$ the degree
that $x$ is {\it totally below} $a$.

\begin{prop} \label{totally below relation is interpolative} {\rm
(Also in \cite{Stubbe05})} Suppose $A$ is a completely
distributive $\Om$-lattice. Then the totally below relation on $A$
is interpolative  in the sense that for all $x, y\in A,$
$$\Downarrow\!(x)(y)= \bv_{z\in
A}\Big(\Downarrow\!(x)(z)*\Downarrow\!(z)(y)\Big).$$
\end{prop}

\noindent{\bf Proof.} \ Let $\lam(w) = \bv_{z\in
A}\Big(\Downarrow\!(x)(z)*\Downarrow\!(z)(w)\Big)$ for all $w\in
A$. Then $\lam\in[A^{\rm op}, \Om]$. Because \begin{eqnarray*}
\sup\lam &=&
\sup\bv_{z\in A}\Big(\Downarrow\!(x)(z)*\Downarrow\!(z)\Big)\\
&=&\sup\bv_{w\in A}\Big(\Downarrow\!(x)(z)*\Downarrow\!(z)(w)\Big)
\otimes w\\
&=& \bv_{w\in A}\bv_{z\in
A}\Big((\Downarrow\!(x)(z)*\Downarrow\!(z)(w))\otimes w\Big)\\
&=& \bv_{z\in A}\Big[\Downarrow\!(x)(z)\otimes\Big({\bv}_{w\in
A}(\Downarrow\!(z)(w)\otimes w\Big)\Big]\\
&=& \bv_{z\in A}\Big[\Downarrow\!(x)(z) \otimes
\sup(\Downarrow\!(z))\Big]\\
&=& \bv_{z\in A}\Big(\Downarrow\!(x)(z)\otimes z\Big)\\
&=& \sup(\Downarrow\!(x)) = x.
\end{eqnarray*} Therefore $\lam\geq\Downarrow\!(x)$ and particularly,
$$\Downarrow\!(x)(y) \leq \bv_{z\in
A}\Big(\Downarrow\!(x)(z)*\Downarrow\!(z)(y)\Big).$$

Conversely, since $\Downarrow: A\lra[A^{\rm op}, \Om]$ is an
$\Om$-functor, for each $z\in A$, we have that
$$\Downarrow\!(x)(z)\leq A(z, x)\leq \bw_{w\in
A}\Big(\Downarrow\!(z)(w) \ra \Downarrow\!(x)(w)\Big) \leq
\Downarrow\!(z)(y) \ra \Downarrow\!(x)(y).$$ Consequently,
$$\Downarrow\!(x)(y)\geq \Downarrow\!(x)(z)*\Downarrow\!(z)(y)$$
for all $z\in A$, and thus,
$$\Downarrow\!(x)(y)\geq \bv_{z\in
A}\Big(\Downarrow\!(x)(z)*\Downarrow\!(z)(y)\Big).$$ $\Box$

The category of all the completely distributive $\Om$-lattices and
the complete $\Om$-lattice morphisms is denoted $\Om$-{\bf CD},
which is a full subcategory of $\Om$-{\bf CLat}.

\begin{eg}Both the singleton discrete $\Om$-category $\bf 1$ and
the terminal $\Om$-category $\top$ are completely distributive
$\Om$-lattices.\end{eg}

\begin{prop}\label{lower sets is CD} {\rm
(Also in \cite{Stubbe05})} Suppose $A$ is an $\Om$-category. Then
$[A^{\rm op},\Om]$ is a completely distributive $\Om$-lattice.
Particularly, both $[\Om^A]= [|A|, \Om]$ and $\Om\cong[\Om^{\bf
1}]$ are completely distributive $\Om$-lattices.
\end{prop}

\noindent{\bf Proof.} \ Let ${\bf y}_A:A\lra[A^{\rm op},\Om]$ be
the Yoneda embedding. By \ref{kan extension}, the $\Om$-functor
${\bf y}_A^\la:[[A^{\rm op},\Om]^{\rm op},\Om]\lra[A^{\rm
op},\Om]$ has a left adjoint. Thus, it suffices to show that the
$\Om$-functors $\sup:[[A^{\rm op},\Om]^{\rm op},\Om]\lra[A^{\rm
op},\Om]$ and ${\bf y}_A^\la$ coincide with each other, i.e.
$\sup\Phi={\bf y}_A^\la(\Phi)=\Phi\circ{\bf y}_A$ for all
$\Om$-functor $\Phi:[A^{\rm op},\Om]^{\rm op}\lra\Om$. Indeed, for
any $x\in A$,
$$\sup\Phi(x)=\bv_{\phi\in[A^{\rm op},\Om]}\Phi(\phi)*\phi(x)
\geq\Phi({\bf
y}_A(x))*{\bf y}_A(x)(x)\geq\Phi({\bf y}_A(x)). $$

On the other hand, because $\Phi$ is an $\Om$-functor from
$[A^{\rm op},\Om]^{\rm op}$ to $\Om$,
$$\sup\Phi(x)=\bv_{\phi\in[A^{\rm op},\Om]}\Phi(\phi)*\phi(x)=
\bv_{\phi\in[A^{\rm op},\Om]}\Phi(\phi)*[A^{\rm op},\Om]({\bf
y}_A(x),\phi)\leq\Phi({\bf y}_A(x)).$$

Therefore,  $\sup(\Phi)=\Phi\circ{\bf y}_A$.\ \ \ $\Box$

\begin{thm}\label{prod is CD} Suppose $\{A_i:i\in J\}$ is a family
of completely distributive $\Om$-lattices. The product
$A=\prod_{i\in J}A_i$ is also a completely distributive
$\Om$-lattice.\end{thm}

\noindent{\bf Proof.}\ \  Step 1. We show that for all
$\phi\in[A^{\rm op},\Om]$, and $j\in J$, the image of $\phi$ under
the projection $p_j:\prod_{i\in J}A_i\lra A_j$ is given by
$p_j(\phi)(t)=\phi(q_j(t))$ for all $t\in A_j$, where
$q_j:A_j\cong A_j\times\prod_{i\neq j}\{0_{A_i}\}\lra\prod_{i\in
J}A_i$ is the left adjoint of the projection $p_j:\prod_{i\in
J}A_i\lra A_j$. In fact, since $\phi$ is  a decreasing function
from the complete lattice $A_0$ to the complete lattice $\Om$, for
all $t\in A_j$ and $x\in p^{-1}(t)$, we have that
$\phi(q_j(t))\geq\phi(x)$. Therefore, $$p_j(\phi)(t)=\bv_{x\in
p_j^{-1}(t)}\phi(x)=\phi(q_j(t)).$$

Step 2. We construct a left adjoint of the $\Om$-functor
$\sup:[A^{\rm op},\Om]\lra A$. For each $j\in J$, let
$d_j:A\lra\Om$ be given by
$$d_j(a)(x)=\left\{\begin{array}{ll}
\Da_j(a_j)(t),   &\mbox{$x=q_j(t)$ for some $t\in A_j$;}\\
0,  &\mbox{otherwise. }
\end{array}
\right.$$

Let $\Da(a)= \ \da\!(\bv_{j\in J}d_j(a))$. Then we claim that
$\Da:A\lra[A^{\rm op},\Om]$ is a left adjoint of $\sup:[A^{\rm
op},\Om]\lra A$. Actually, for all $a\in A$ and $\phi\in[A^{\rm
op},\Om]$,
\begin{eqnarray*}
[A^{\rm op},\Om](\Da(a),\phi) &=&[A^{\rm
op},\Om]\Big(\da\!\Big({\bv}_{j\in J}d_j(a)\Big),\phi\Big)\\
&=&[\Om^A]\Big({\bv}_{j\in J}d_j(a),\phi\Big)\ \ \
{\rm (Example\ \ref{lower is left adjoint})}\\
&=&\bw_{j\in J}\bw_{x\in A}d_j(a)(x)\ra\phi(x)\\
&=&\bw_{j\in J}\bw_{t\in A_j}\Da_j(a_j)(t)\ra p_j(\phi)(t)\\
&=&\bw_{j\in J}[A_j^{\rm op},\Om](\Da_j(a_j),p_j(\phi))\\
&=&\bw_{j\in J}A_j(a_j,\sup{_{A_j}}p_j(\phi))\\
&=&A(a,\sup\phi),\end{eqnarray*} where the last equality holds
because the complete $\Om$-lattice morphism $p_j$ preserves sups,
i.e. $\sup_{A_j}p_j(\phi)=p_j(\sup_A\phi).$ Therefor, $\Da$ is a
left adjoint of $\sup:[A^{\rm op},\Om]\lra A$ as claimed. \ \
$\Box$

\section{Subalgebras and quotient algebras}

\begin{df} Suppose that $A$ is an $\Om$-category and $B$ is a
subcategory of $A$. Then $B$ is said to be a subalgebra of $A$ if
the embedding functor $i:B\lra A$ preserves both sups and
infs.\end{df}

Suppose $A$ is a complete $\Om$-lattice and $B$ is a subcategory
of $A$. It is routine to check that the following conditions are
equivalent: (1) $B$ is a subcategory of $A$; (2) $B$ is closed
with respect to tensors and cotensors in $A$ and $B_0$ is closed
with respect to joins and meets in $A_0$;  (3) the embedding
functor $i:B\lra A$ has both a left and a right adjoint.

By (2), every subalgebra $B$ of a complete $\Om$-lattice $A$ is
itself a complete $\Om$-lattice. Here is another proof of this
fact. Let $i:B\lra A$ denote the embedding and $k:A\lra B$ be a
left adjoint of $i$. Then the Yoneda embedding $\y_B:B\lra [B^{\rm
op},\Om]$ can be written as a composition $i^\la\circ\y_A\circ i:
B\lra A\lra[A^{\rm op}, \Om]\lra[B^{\rm op}, \Om]$. Thus, $\y_B$
has a left adjoint given by $\sup_B= k\circ\sup_A\circ i^\ra_\ell:
[B^{\rm op}, \Om]\lra[A^{\rm op}, \Om]\lra A\lra B$.

Moreover, for completely distributive $\Om$-lattices, we have the
following.

\begin{prop}\label{sub lat is CD} Suppose $L$ is a completely
distributive $\Om$-lattice and $M$ is a subalgebra of $L$. Then
$M$ is also completely distributive.
\end{prop}

\noindent{\bf Proof.}\ \  Suppose that $i:M\lra L$ is the
corresponding embedding and  $k:L\lra M$ is a left adjoint of $i$.
Let $\Da_L:L\lra[L^{\rm op},\Om]$ be a left adjoint of
$\sup_L:[L^{\rm op},\Om]\lra L$. Then, we say that $\Da_M =
k^\ra_\ell\circ\Da_L\circ i: M\lra L\lra[L^{\rm op},
\Om]\lra[M^{\rm op}, \Om]$ is a left adjoint of $\sup_M:[M^{\rm
op},\Om]\lra \Om$, hence, $M$ is completely distributive.

Indeed, for all $x\in A$ and $\phi\in [M^{\rm op},\Om]$,
$$[M^{\rm op},\Om](\Da_M(x),\phi]=[M^{\rm op},\Om](k^\ra_\ell
\circ\Da_L\circ i(x),\phi)=L(i(x),\sup{_L}\circ k^\la(\phi)).$$
And
$$M(x,\sup{_M}(\phi))=L(i(x), i(\sup{_M} \phi ))
=L(i(x),\sup{_L}\circ i^\ra_\ell(\phi))$$ since $i:M\lra L$
preserves sups.

So, it suffices to show that $k^\la(\phi)=i^\ra_\ell(\phi)$. In
fact, for all $y\in A$,
\begin{eqnarray*}i^\ra_\ell(\phi)(y)&=&\bv_{x\in
M}\phi(x)*L(y,i(x))=\bv_{x\in M}\phi(x)*M(k(y), x)\\
&\geq&\phi(k(y))*M(k(y),k(y))\geq\phi(k(y))=k^\la(\phi)(y).
\end{eqnarray*}
On the other hand, since $\phi:M^{\rm op}\lra \Om$ is an
$\Om$-functor, \begin{eqnarray*}i^\ra_\ell(\phi)(y)&=&\bv_{x\in
M}\phi(x)*L(y,i(x))=\bv_{x\in M}\phi(x)*M(k(y), x)\\
&\leq&\phi(k(y))=k^\la(\phi)(y).
\end{eqnarray*}
Thus, $k^\la(\phi)=i^\ra_\ell(\phi)$ as desired. \ \ $\Box$

\begin{cor} $\Om$-{\bf CD} is a complete category.\end{cor}

\noindent{\bf Proof.} \  This is because that $\Om$-{\bf CD} has a
terminal object and products by 4.4 and \ref{prod is CD}, and that
it has equalizers by \ref{sub lat is CD} and \ref{equalizer
exist}.\ \ \ $\Box$

In the following we shall show that the subalgebras of a
completely distributive $\Om$-lattice $A$ can be equivalently
described by certain closure operators on $A$.

\begin{df}A closure operator on an $\Om$-category $A$ is an
idempotent $\Om$-functor $c:A\lra A$ such that  $A(x, c(x))\geq I$
for all $x\in A$.\end{df}

\begin{lem} Suppose $A$ is a complete $\Om$-lattice and
$B\subseteq A$ is a subalgebra of $A$. Let $c: A\lra A$ be given
by $c(x) = \bw\{y\in B\ |\ x\leq y\}$. Then $c$ is a cocontinuous
closure operator on $A$.\end{lem}

\noindent{\bf Proof.} \ (1)  $c$  is an $\Om$-functor. It suffices
to show that $\alpha\otimes c(x)\leq c(\alpha\otimes x)$ for all
$\alpha\in\Om$ and $x\in A$. Indeed, for all $z\in B$, if
$\alpha\otimes x\leq z$, then $x\leq\alpha\rightarrowtail z$ by
\ref{property of cotensor}(1). Since $B$ is a subalgebra of $A$,
$\alpha\rightarrowtail z\in B$. Hence, $c(x)\leq
\alpha\rightarrowtail z$ and  $\alpha\otimes c(x)\leq
c(\alpha\otimes x)$ by arbitrariness of $z$.

(2) For every $x\in A$, $A(x, c(x))\geq I$. This is trivial by
definition.

(3) $c: A_0\lra A_0$ preserves joins. Suppose that $(x_t)_{t\in
T}\subseteq A$. Then
$$c\Big({\bv}_{t\in T}x_t\Big)= c\circ c\Big({\bv}_{t\in
T}x_t\Big)\geq c\Big({\bv}_{t\in T}c(x_t)\Big)\geq {\bv}_{t\in
T}c(x_t).$$ On the other hand, since $\bv_{t\in T}x_t\leq
\bv_{t\in T}c(x_t)$ and $\bv_{t\in T}c(x_t)\in B$ by assumption,
we have that $c\Big(\bv_{t\in T}x_t\Big)\leq \bv_{t\in T}c(x_t)$.
Hence $c\Big(\bv_{t\in T}x_t\Big)= \bv_{t\in T}c(x_t)$.

(4) $c$ preserves tensors, i.e., $c(\alpha\otimes x)
=\alpha\otimes c(x)$ for all $\alpha\in\Om$ and $x\in A$. By (1),
we need only check that $c(\alpha\otimes x) \leq\alpha\otimes
c(x)$. By definition, $c(x)$ is the least element in $B$ such that
$x\leq c(x)$. Thus, $\alpha\otimes x\leq \alpha\otimes c(x)$.
Consequently, $c(\alpha\otimes x) \leq\alpha\otimes c(x)$ since
$\alpha\otimes c(x)\in B$.

Therefore, $c$ is a cocontinuous closure operator. \ \ \ $\Box$

\begin{lem}Suppose that $A$ is a complete $\Om$-lattice and
$c: A\lra A$ is a cocontinuous closure operator. Then $c(A) =
\{c(x)\ |\ x\in A\}$ is a subalgebra of $A$.\end{lem}

\noindent{\bf Proof.} \ (1) Suppose $(x_t)_{t\in T}\subseteq
c(A)$. Then $$c\Big({\bv}_{t\in T}x_t\Big) = {\bv}_{t\in T}c(x_t)
= {\bv}_{t\in T}x_t$$ and
$${\bw}_{t\in T}x_t \leq c\Big({\bw}_{t\in T}x_t\Big) \leq {\bw}_{t\in
T}c(x_t) = {\bw}_{t\in T}x_t. $$ Thus, $(c(A))_0$ is a complete
sublattice of $A_0$.

(2) For all $\alpha\in\Om$ and $x\in c(A)$, $\alpha\otimes x\in
c(A)$. This is because $$\alpha\otimes x\leq c(\alpha\otimes x)
=\alpha\otimes c(x) =\alpha\otimes x.$$

(3) For all $\alpha\in\Om$ and $x\in c(A)$, $\alpha\rightarrowtail
x\in c(A)$, i.e., $c(\alpha\rightarrowtail x) = \alpha
\rightarrowtail x$. Since $$A(c(y), x)\leq A(y, x)\leq A(c(y),
c(x)) = A(c(y), x)$$ for all $y\in A$, we have that
\begin{eqnarray*} \alpha\rightarrowtail x &=& \bv\{y\in A \ |\
\alpha\leq A(y,
x)\}\\
&=& \bv\{y\in A \ |\ \alpha\leq A(c(y),x)\}\\
&=& \bv\{c(y) \ |\ y\in A\ {\rm and}\ \alpha\leq A(y,x)\}\\
&=&c(\alpha\rightarrowtail x).\end{eqnarray*}

Therefore, $c(A)$ is a subalgebra of $A$. \ \ \ $\Box$

A combination of the above two lemmas yields the following.

\begin{prop}The subalgebras of a completely distributive
$\Om$-lattice $A$ correspond bijectively to the cocontinuous
closure operators on $A$. \ \ \ $\Box$
\end{prop}

\begin{df} Let $A$ and $B$ be $\Om$-categories. $B$ is
said to be a quotient algebra of $A$ if there is a surjective
$\Om$-functor $f:A\lra B$ such that $f$ has both a left adjoint
and a right adjoint.
\end{df}

Suppose that $A$ is a complete $\Om$-lattice and $B$ is a quotient
algebra of $A$. By definition, there is a surjective functor
$f:A\lra B$ which has both a left adjoint and a right adjoint. Let
$k:L\lra M$ be a right adjoint of $f$. Then $k$ is an isometry and
the Yoneda embedding $\y_B:B\lra [B^{\rm op},\Om]$ can be written
as a composition $k^\la\circ\y_A\circ k: B\lra A\lra[A^{\rm op},
\Om]\lra[B^{\rm op}, \Om]$. Thus, $\y_B$ has a left adjoint given
by $\sup_B= f\circ\sup_A\circ f^\la: [B^{\rm op}, \Om]\lra[A^{\rm
op}, \Om]\lra A\lra B$. Therefore, $B$ is a complete
$\Om$-lattice. To conclude, every quotient algebra of a complete
$\Om$-lattice is also a complete $\Om$-lattice.

\begin{prop}\label{image is CD} Every quotient
algebra of a completely distributive $\Om$-lattice is also a
completely distributive $\Om$-lattice.\end{prop}

\noindent{\bf Proof.} \  Suppose that $L$ is a completely
distributive $\Om$-lattice  and $f:L\lra M$ is a surjective
$\Om$-functor  with a left adjoint $j: M\lra L$ and a right
adjoint $k: M\lra L$.  Then, $M$ is a complete $\Om$-lattice by
the above observation and $j,k$ are injective isometric
$\Om$-functors by \ref{injec is embed}. Let
$\Da_M=f^\ra_\ell\circ\Da_L\circ j$, where $\Da_L$ is the left
adjoint of $\sup{_L}$. We left it to the reader to check that
$\Da_M$ is a left adjoint of $\sup_M$. Hence, $M$ is completely
distributive. \ \ $\Box$

The following classical characterization of complete
distributivity was established by Raney and B\"{u}chi
independently in \cite{Buchi, Raney52, Raney53}.

\begin{thm}\label{RB}{\rm(Raney-B\"{u}chi)}
A complete $\Om$-lattice $L$ is completely distributive if and
only if there is a set $X$ such that $L$ is a quotient algebra of
some subalgebra of $[\Om^X]$.\end{thm}

\noindent{\bf Proof.}\ \  Sufficiency: This follows from
\ref{lower sets is CD}, \ref{sub lat is CD} and \ref{image is CD}
immediately.

Necessity: since $[L^{\rm op},\Om]$ is a subcategory of $[\Om^L]$
and the inclusion $i:[L^{\rm op},\Om]\lra[\Om^L]$ has both a left
adjoint and a right adjoint, $[L^{\rm op},\Om]$ is a subalgebra of
$[\Om^L]$. Since $\sup: [L^{\rm op},\Om]\lra L$ is surjective and
has both a left adjoint $\Da: L\lra[L^{\rm op},\Om]$ and a right
adjoint $\y: L\lra[L^{\rm op},\Om]$, $L$ is a quotient algebra of
$[L^{\rm op},\Om]$. \ \ $\Box$

\begin{df}A kernel operator on an $\Om$-category $A$ is an
idempotent $\Om$-functor $k:A\lra A$ such that  $A(k(x), x)\geq I$
for all $x\in A$.\end{df}

Let $A$ be a complete $\Om$-lattice and $k: A\lra A$ be a
cocontinuous kernel operator on $A$. Define $x\sim y$ in $A$ if
$k(x)=k(y)$. Then

(1) If $x_t\sim y_t$ for all $t\in T$, then $\bv_{t\in
T}x_t\sim\bv_{t\in T}y_t$. This is trivial since $k$ preserves
joins.

(2) If $x\sim y$, then $\alpha\otimes x\sim \alpha\otimes y$ for
all $\alpha\in\Om$ because $k$ preserves tensors.

(3) If $x_t\sim y_t$ for all $t\in T$, then $\bw_{t\in
T}x_t\sim\bw_{t\in T}y_t$. To see this, note at first that
$$k\Big({\bw}_{t\in T}x_t\Big)\leq {\bw}_{t\in T}k(x_t)=
{\bw}_{t\in T}k(y_t)\leq {\bw}_{t\in T}y_t.$$ So
$$k\Big({\bw}_{t\in T}x_t\Big)= k\circ k\Big({\bw}_{t\in
T}x_t\Big)\leq k\Big({\bw}_{t\in T}y_t\Big).$$ Exchanging the
roles of $x$ and $y$ we obtain that
$$k\Big({\bw}_{t\in T}y_t\Big) \leq k\Big({\bw}_{t\in
T}x_t\Big).$$ Therefore, $$k\Big({\bw}_{t\in T}y_t\Big) =
k\Big({\bw}_{t\in T}x_t\Big).$$

(4) If $x\sim y$, then, for all $\alpha\in\Om$,
$\alpha\rightarrowtail x\sim \alpha\rightarrowtail y$. That is,
$k(\alpha\rightarrowtail x)=k(\alpha\rightarrowtail y)$.

To see this, we assert at first that for all $\alpha\in\Om$ and
$x\in A$, the set $\{k(z)\in A \ |\ \alpha\leq A(z, x)\}$ equals
the set $\{k(z)\in A \ |\ \alpha\leq A(k(z), k(x))\}$.

The inclusion $\{k(z)\in A \ |\ \alpha\leq A(z, x)\}\subseteq
\{k(z)\in A \ |\ \alpha\leq A(k(z), k(x))\}$ is trivial since $k$
is an $\Om$-functor. For the converse inclusion, suppose that
$\alpha\leq A(k(z), k(x))$. Then $\alpha\leq A(k(z), k(x))*A(k(x),
x)\leq A(k(z), x)$. Thus, $k(z) = k(k(z))\in \{k(z)\in A \ |\
\alpha\leq A(z, x)\}.$

Therefore,
\begin{eqnarray*} k(\alpha\rightarrowtail
x) &=& k\Big(\bv\{z\in A \ |\ \alpha\leq A(z, x)\}\Big)\\
&=& \bv\{k(z)\in A \ |\ \alpha\leq A(z, x)\}\\
&=& \bv\{k(z)\in A \ |\ \alpha\leq A(k(z), k(x))\}\\
&=& \bv\{k(z)\in A \ |\ \alpha\leq A(k(z), k(y))\}\\
&=& \bv\{k(z)\in A \ |\ \alpha\leq A(z, y)\}\\
&=& k\Big(\bv\{z\in A \ |\ \alpha\leq A(z, y)\}\Big)\\
&=& k(\alpha\rightarrowtail y).
\end{eqnarray*}

By (1)-(4) in the above we see that the subset $R =\{(x, y)\ |\
x\sim y\}\subseteq  A\times A$ is not only an equivalence relation
on $A$, but also a subalgebra of $A\times A$. Let $B = A/\sim$ and
$q: A\lra B$ be the corresponding quotient map. We define an
$\Om$-category structure on $B$ as follows: for all $[x], [y]\in
B$, let $B([x], [y]) = A(k(x), k(y))$. Then it is easy to verify
that $B$ becomes an $\Om$-category and $q:A\lra B$ is an
$\Om$-functor. An we leave it to the reader to check that $B$ is a
complete $\Om$-lattice and that $q$ preserves joins, tensors,
meets, and cotensors. Thus, $q$ has, at the same time, a left
adjoint and a right adjoint. That means, $B$ is a quotient algebra
of $A$. Therefore, every cocontinuous kernel operator on a
complete $\Om$-lattice $A$ is associated with a quotient algebra
of $A$.

On the other hand, suppose that $B$ is a quotient algebra of a
complete $\Om$-lattice $A$ with the quotient map $q: A\lra B$. By
definition, $q$ has a left adjoint $f: B\lra A$. Let $k=f\circ q:
A\lra B\lra A$. Then $k$ is a cocontinuous kernel operator on $A$.

We leave it to the reader to check that the above processes from
quotient algebras of a complete $\Om$-lattice $A$ to the
cocontinuous kernel operators on $A$ and vice versa are inverse to
each other. Particularly, we arrive at the following.

\begin{prop}The quotient algebras of a completely distributive
$\Om$-lattice $A$ correspond bijectively to the cocontinuous
kernel operators on $A$. \ \ \ $\Box$
\end{prop}

For each complete $\Om$-lattice $A$, $A_0$ must be a complete
lattice. But, the complete distributivity of $A$ (as an
$\Om$-lattice) does not imply the complete distributivity of
$A_0$. However, we have the following.

\begin{prop} The following two conditions are equivalent:

{\rm (1)} The complete lattice $\Om$ is completely distributive.

{\rm (2)} For every completely distributive $\Om$-lattice $A$,
$A_0$ is completely distributive.
\end{prop}

\noindent{\bf Proof.}\ \  $(2)\Rightarrow(1)$: Since $(\Om,\ra)$
is a completely distributive $\Om$-lattice, thus, $\Om= (\Om,
\ra)_0$ is a completely distributive lattice.

$(1)\Rightarrow(2)$: Suppose $\Om$ is a completely distributive
complete lattice. Then, $[\Om^A]_0$ is a completely distributive
complete lattice. Since $\sup: [A^{\rm op},\Om]\lra A$ preserves
sups and infs, $A$ is a quotient algebra of $[A^{\rm op},\Om]$,
which is a subalgebra of $[\Om^A]$. Therefore, $A_0$ is a quotient
algebra of  a subalgebra $[A^{\rm op},\Om]_0$ of $[\Om^X]_0$.
Letting $\Om={\bf 2}$ in \ref{RB}, we obtain that $A_0$ is a
completely distributive lattice.\ \ $\Box$

\section{The $\Om$-category of left adjoints}

In this section, we discuss the complete distributivity of functor
categories. The main result is that if $A$ and $B$ are completely
distributive $\Om$-lattices, then so is the subcategory of the
functor category $[A, B]$ consisting of left adjoints from $A$ to
$B$.

Suppose $A$ and $B$ are  $\Om$-categories. Let $[A\ra_\ell B]$
denote the subcategory of $[A, B]$ consisting of left adjoint
functors, i.e., every element in $[A\ra_\ell B]$ has a right
adjoint.

\begin{prop} Suppose $A$ and $B$ are complete $\Om$-lattices.
Then $[A\ra_\ell B]$ is a complete $\Om$-lattice.\end{prop}

\noindent{\bf Proof.} \ (1) Suppose that $(f_t)_{t\in
T}\subseteq[A\ra_\ell B]$, the pointwise join $f = \bv_{t\in
T}f_t$ given by $f(x)=\bv_{t\in T}f_t(x)$ is a cocontinuous
$\Om$-functor $A\lra B$, where the join is taken in the complete
lattice $B_0$. This is because that (i) $$f(\alpha\otimes x) =
\bv_{t\in T}f_t(\alpha\otimes x)= \bv_{t\in T}\alpha\otimes
f_t(x)= \alpha\otimes\bv_{t\in T}f_t(x)=\alpha\otimes f(x),$$
i.e., $f$ preserves tensors; and (ii) $$f({\bv}_{i\in J}x_i) =
\bv_{t\in T}f_t({\bv}_{i\in J}x_i)= \bv_{i\in J}\bv_{t\in T}
f_t(x_i)= \bv_{i\in J}f(x_i),$$ i.e., $f: A_0\lra B_0$ preserves
joins. So, $[A\ra_\ell B]_0$ is a complete lattice.

(2) For every $f\in[A\ra_\ell B]$ and $\alpha\in\Om$, the tensor
$\alpha\otimes f$ of $\alpha$ and $f$ in $[A, B]$ is cocontinuous,
i.e., $[A\ra_\ell B]$ is closed under tensor in $[A, B]$. The
proof is trivial and hence omitted here.

Therefore, $[A\ra_\ell B]$ is a complete $\Om$-lattice with the
tensor inherited from $[A, B]$. \ \ \ $\Box$

It should be noted that though $[A\ra_\ell B]$ is a complete
$\Om$-lattice, it is not necessarily a subalgebra of $[A, B]$.

\begin{lem}\label{lemma} Suppose $A$ is an $\Om$-category
and $B$ is a completely distributive $\Om$-lattice. Then $[A,B]$
is also a completely distributive $\Om$-lattice.\end{lem}

\noindent{\bf Proof.} \ Indeed, $[A,B]$ is a subalgebra of the
completely distributive $\Om$-lattice $[|A|,B]=B^{|A|}$. To this
end, we need only check that (1) for all $\alpha\in\Om$ and
$f\in[A, B]$, the tensor $\alpha\otimes f$ and cotensor
$\alpha\rightarrowtail f$ in $B^{|A|}$ are $\Om$-functors; and (2)
for every family $(f_t)_{t\in T}\subseteq[A, B]$, the pointwise
join $\bv_{t\in T}f_t$ and pointwise meet $\bw_{t\in T}f_t$ are
$\Om$-functors.  The details are left to the reader. \ \ \ $\Box$

\begin{thm} Suppose $A$ and $B$ are completely distributive
$\Om$-lattices.  Then $\Om$-category $[A\ra_\ell B]$ of left
adjoints is a completely distributive $\Om$-lattice.\end{thm}

\noindent{\bf Proof.} \ Our strategy is to show that $[A\ra_\ell
B]$ is a quotient algebra of $[A, B]$, which is a completely
distributive $\Om$-lattice by Lemma \ref{lemma}.

Define $k: [A, B]\lra [A, B]$ by $k(f)(a) = \bv_{x\in
A}\Da\!(a)(x)\otimes f(x)$ for all $f\in[A, B]$ and $a\in A$,
where, $\otimes$ denotes the tensor in $B$.

(1) $k: [A, B]_0\lra[A, B]_0$ preserves order. Trivial by
definition.

(2) $k(f)\leq f$ for all $f\in[A, B]$. Because $\Da\!(a)(x)\leq
A(x, a)$ for $a, x\in A$, then $\Da\!(a)(x)\otimes f(x)\leq A(x,
a)\otimes f(x)\leq f(a).$

(3) $k$ is idempotent. For all $f\in[A, B]$ and $a\in A$,
\begin{eqnarray*} k\circ k(f)(a) &=& \bv_{x\in A}\Da\!(a)(x)\otimes
k(f)(x)\\
&=& \bv_{x\in A}\Big[\Da\!(a)(x)\otimes \Big(\bv_{y\in
A}\Da\!(x)(y)\otimes f(y)\Big)\Big]\\
&=& \bv_{y\in A}\Big[\Big(\bv_{x\in
A}\Da\!(a)(x)*\Da\!(x)(y)\Big)\otimes f(y)\Big]\\
&=& \bv_{y\in A}\Da\!(a)(y)\otimes f(y)\ \ \
({\rm Prop. }\ \ref{totally below relation is interpolative})\\
&=& k(f)(a).
\end{eqnarray*}

(4) $k$ preserves tensors. \begin{eqnarray*}k(\alpha\otimes f)(a)
&=& \bv_{x\in A}[\Da\!(a)(x)\otimes (\alpha\otimes f)(x)]\\
&=& \bv_{x\in A}[\Da\!(a)(x)\otimes (\alpha\otimes f(x))]\\
&=& \bv_{x\in A}[\alpha\otimes(\Da\!(a)(x)\otimes  f(x))]\\
&=& \alpha\otimes\bv_{x\in A}[\Da\!(a)(x)\otimes  f(x)]\\
&=& \alpha\otimes k(f)(a)\\
&=& (\alpha\otimes k(f))(a).
\end{eqnarray*}

(5) $k: [A, B]_0\lra[A, B]_0$ preserves joins. Indeed, for all
$a\in A$,
\begin{eqnarray*}k\Big(\bv_{t\in T}f_t\Big)(a)
&=& \bv_{x\in A}\Big[\Da\!(a)(x)\otimes\Big(
\bv_{t\in T}f_t\Big)(x)\Big]\\
&=& \bv_{x\in A}\Big[\Da\!(a)(x)\otimes
\bv_{t\in T}f_t(x)\Big]\\
&=& \bv_{t\in T}\bv_{x\in A}\Big[\Da\!(a)(x)\otimes f_t(x)\Big]\\
&=& \bv_{t\in T}k(f_t)(a).
\end{eqnarray*}

(6) $k$ is an $\Om$-functor. This follows from a combination of
(1) and (4).

(7) $k(f) = f$ if and only if $f$ is cocontinuous. Necessity is
trivial by (4)-(6).

Sufficiency: If $f$ is cocontinuous, then for all $a\in A$,
\begin{eqnarray*}f(a) &=& f(\sup\Da\!(a))=\sup f(\Da\!(a))\\
&=& \bv_{y\in B}f(\Da\!(a))(y)\otimes y \ \ \ {\rm
(Prop. \ \ref{sup by join and tensor}(1))}\\
&=& \bv_{y\in B}\Big[\Big(\bv_{f(x)=y}\Da\!(a)(x)\Big)\otimes y\Big]\\
&=& \bv_{x\in A}\Da\!(a)(x)\otimes f(x)\\
&=& k(f)(a).
\end{eqnarray*}

Therefore, $k$ is a cocontinuous kernel operator on $[A, B]$, the
corresponding quotient algebra is $[A\ra_\ell B]$ by (7). Hence
$[A\ra_\ell B]$ is completely distributive. \ \ \ $\Box$

A natural question is to ask whether the category of right
adjoints between two completely distributive $\Om$-lattices is
also completely distributive. At first, we say that for any
$\Om$-categories $A$ and $B$, the $\Om$-category $[B\ra_r A]$ of
right adjoints from $B$ to $A$ (as a subcategory of $[B, A]$) is
isomorphic to the dual category of $[A \ra_\ell B]$ consisting of
left adjoints from $A$ to $B$. To see this, we need only check
that for any $\Om$-adjunctions $(f_1, g_1)$ and $(f_2, g_2)$,
$$\bw_{x\in A}B(f_1(x), f_2(x)) = \bw_{y\in B}A(g_2(y), g_1(y)).$$
Indeed, for any $y\in B$, let $x=g_2(y)$. Then
\begin{eqnarray*}B(f_1(x), f_2(x))&=& B(f_1(x), f_2(g_2(y)))
\leq B(f_1(x),y)\\
& =&  A(x, g_1(y)) = A(g_2(y), g_1(y)).\end{eqnarray*} Therefore,
$$\bw_{x\in A}B(f_1(x), f_2(x)) \leq \bw_{y\in B}A(g_2(y), g_1(y)).$$

Conversely, for any $x\in A$, let $y=f_2(x)$. Then
\begin{eqnarray*}A(g_2(y), g_1(y))&=& A(g_2(f_2(x)), g_1(y))
\leq A(x,g_1(y))\\
& =&  B(f_1(x), y) = B(f_1(x), f_2(x)).\end{eqnarray*} Therefore,
$$\bw_{x\in A}B(f_1(x), f_2(x)) \geq \bw_{y\in B}A(g_2(y), g_1(y)).$$

Particularly, if both $A$ and $B$ completely distributive
$\Om$-lattices, then the dual category of the $\Om$-category of
right adjoints from $B$ to $A$, or that from $A$ to $B$, is a
completely distributive $\Om$-lattice. But, this does not mean
that $[B\ra_r A]$ is completely distributive since, as we shall
see in the next section, the dual category of a completely
distributive $\Om$-lattice is not necessarily a completely
distributive $\Om$-lattice.

\section{When $\Om$ is a Girard quantale}

In this section, we investigate the complete distributivity of the
dual of a completely distributive $\Om$-lattice. The result shows
that this depends heavily on the properties of $\Om$. That is,
when $\Om$ is an integral commutative quantale, then every
completely distributive $\Om$-lattice is dually completely
distributive if and only if $\Om$ is a Girard quantale. This
conclusion should be compared with the fact that in any topos
$\cal E$, the dual of every constructive completely distributive
lattice is constructive completely distributive if and only if
$\cal E$ is a Boolean topos \cite{RW2,W}.

By definition, it is easy to see that an $\Om$-category $A$ is
dually completely distributive, i.e., $A^{\rm op}$ is completely
distributive, if and only if the $\Om$-functor $\inf:[A,\Om]^{\rm
op}\lra A$ has a right adjoint.

\begin{eg}Suppose that $\Om$ is a commutative Girard quantale. Then
the $\Om$-functor $\inf:[\Om,\Om]^{\rm op}\lra\Om$ has a right
adjoint. Thus $\Om^{\rm op}$ is a completely distributive
$\Om$-lattice. To this end, we show that the mapping
$d:\Om\lra[\Om,\Om]^{\rm op}$, given by $d(x)(t)=x\ra 0$ for all
$x\in\Om,\ t\in \Om$, is a right adjoint of $\inf:[\Om,\Om]^{\rm
op}\lra\Om$.

At first, for all $\psi\in [\Om,\Om]$, we have that

(1) $\psi$ is  increasing; and

(2) for any $t\in\Om$, $\psi(0)\ra 0\leq\psi(t)\ra t$.

To see (2), for any $t\in\Om$,\begin{eqnarray*} t\ra
0\leq\psi(t)\ra \psi(0)&\iff&
\psi(t)\leq(t\ra 0)\ra\psi(0)=(\psi(0)\ra 0)\ra t\\
&\iff& \psi(0)\ra 0\leq\psi(t)\ra t.\end{eqnarray*}

Therefore, for all $\psi\in [\Om,\Om]$,
$$\inf\psi=\bw_{t\in\Om}\psi(t)\ra t=\psi(0)\ra 0.$$

Hence, for all $\psi\in[\Om,\Om]$ and  $x\in \Om$,
\begin{eqnarray*}\Om(\inf\psi,x)&=&(\psi(0)\ra 0)\ra x
=(x\ra 0)\ra\psi(0)\\
&=&(x\ra 0)\ra\bw_{t\in \Om}\psi(t)=\bw_{t\in\Om}(x\ra
0)\ra\psi(t)\\
&=& [\Om,\Om]^{\rm op}(\psi,d(x)).
\end{eqnarray*} That means, $d$ is a right adjoint of
$\inf$.\end{eg}

Recall that an $\Om$-subset $\phi: A\lra \Om$ of an $\Om$-category
$A$ is said to be finite if the set $\{x\in A\ |\ \phi(x)\not=0\}$
is finite. Noticing that a complete lattice $H$ is a complete
Heyting algebra if and only if the supremum operator $\sup:{\cal
D}(H)\lra H$ preserves finite meets, we introduce the following.

\begin{df} A complete $\Om$-lattice $A$ is called a complete
$\Om$-Heyting algebra if the supremum operator $\sup:[A^{\rm
op},\Om] \lra A$ preserves finite infs.
\end{df}

Clearly, for a complete $\Om$-lattice $A$, the dual $A^{\rm op}$
is a complete $\Om$-Heyting algebra if and only if the
$\Om$-functor $\inf:[A,\Om]^{\rm op}\lra A$ preserves finite sups.

Every completely distributive $\Om$-lattice is a complete
$\Om$-Heyting algebra because the supremum operator has a left
adjoint, hence it preserves (all) infs. Consequently, the complete
$\Om$-lattices $[A^{\rm op},\Om]$, $[\Om^X]$ and $(\Om,\ra)$ are
all complete $\Om$-Heyting algebras.

\begin{prop}Suppose $A$ is a complete  $\Om$-lattice.Then, the
followings are equivalent:

$(1)$ $A$ is a complete $\Om$-Heyting algebra.

$(2)$ $\sup: [A^{\rm op}, \Om]_0\lra A_0$ preserves finite meets
and $\sup(\alpha\ra\lam) =\alpha\ra\sup\lam$ for all
$\alpha\in\Om, \lam\in[A^{\rm op}, \Om]$.\ \ \ $\Box$
\end{prop}

Now we are at the position to prove the main result of this
section, Theorem 1.1 stated in the introduction.

\noindent{\bf Proof of Theorem 1.1.} $(1)\Rightarrow (2)$: Suppose
$\Om$ is a commutative Girard quantale. We claim at first that for
any complete $\Om$-lattice $L$, there is an isomorphism $[L^{\rm
op},\Om]\cong[L,\Om]^{\rm op}$. In fact, define $\neg:[L^{\rm
op},\Om]\lra[L,\Om]^{\rm op}$ by $\neg\phi(x)=\phi(x)\ra 0$ for
all $\phi\in[L^{\rm op},\Om]$ and $x\in A$. At first, $\neg$ is
bijective by the law of double negation. Secondly, because
$$(\alpha\ra 0)\ra(\beta\ra 0)=\beta\ra((\alpha\ra 0)\ra 0)=
\beta\ra \alpha$$ for all $\alpha,\beta\in \Om$, we obtain that
for all $\phi_1,\phi_2\in [L^{\rm op},\Om]$,
\begin{eqnarray*}[L^{\rm op},\Om](\phi_1,\phi_2)&=&\bw_{x\in
L}\phi_1(x)\ra\phi_2(x)\\
&=&\bw_{x\in L}(\phi_2(x)\ra 0)\ra(\phi_1(x)\ra 0)\\
&=&[L,\Om]^{\rm op}(\phi_1\ra 0,\phi_2\ra 0).\end{eqnarray*}
Therefore, $\neg:[L^{\rm op},\Om] \lra[L,\Om]^{\rm op}$ is an
$\Om$-isomorphism and $[L^{\rm op},\Om]^{\rm op}$ is a completely
distributive $\Om$-lattice since  $[L,\Om]$ is completely
distributive by Proposition \ref{lower sets is CD}.

Since $L$ is completely distributive, $\sup:[L^{\rm op},\Om]\lra
L$ is a complete $\Om$-lattice morphism. Consequently, $\sup^{\rm
op}:[L^{\rm op},\Om]^{\rm op}\lra L^{\rm op}$ is also a complete
$\Om$-lattice morphism. Thus,  as a quotient algebra of a
completely distributive $\Om$-lattice, $L^{\rm op} $ is completely
distributive.

$(2)\Rightarrow (3)$: Trivial.

$(3)\Rightarrow (1)$: Suppose $\Om^{\rm op}$ is a complete
$\Om$-Heyting algebra. Then the $\Om$-functor $\inf:[\Om,\Om]^{\rm
op}\lra\Om$ preserves finite sups. Particularly, $\inf$ preserves
tensors. Therefore, for each $\alpha\in\Om$, $$\inf(\alpha\otimes
\underline{0}) = \alpha*\inf\underline{0} = \alpha * 1=\alpha,$$
where $\alpha\otimes \underline{0}$ denotes the tensor of $\alpha$
and the constant functor $\underline{0}$ in $[\Om,\Om]^{\rm op}$,
or equivalently, the cotensor of $\alpha$ and $\underline{0}$ in
$[\Om,\Om]$, that is, $\alpha\otimes \underline{0} =
\alpha\ra\underline{0}$. Meanwhile, $$\inf(\alpha\ra\underline{0})
= \bw_{x\in A}((\alpha\ra 0)\ra x) = (\alpha\ra 0)\ra 0.$$

Therefore, $\alpha = (\alpha\ra 0)\ra 0$. Hence, $\Om$ is a
commutative Girard quantale.\ \ $\Box$

\begin{cor}Suppose $(H,\wedge,\ra)$ is a complete Heyting algebra.
 Then $(H,\ra)^{\rm op}$ is an $H$-Heyting algebra if and only if
 $H$ is a boolean algebra.\ \ $\Box$
\end{cor}

\begin{cor} Suppose $\Om$ is a
commutative Girard quantale, $A$ and $B$ are completely
distributive $\Om$-lattices. Then the category of right adjoints
$[A\ra_r B]$ between $A$ and $B$ is a completely distributive
$\Om$-lattices.\ \ $\Box$ \end{cor}

We end this section with a conclusion on the existence of free
completely distributive $\Om$-lattices when $\Om$ is a commutative
Girard quantale.

\begin{prop} If $\Om$ is a commutative Girard  quantale, then
the forgetful functor $G: \Om$-{\bf CD}$\lra {\bf Set}$ has a left
adjoint. Hence, for every set $X$, there is free completely
distributive $\Om$-lattice generated by $X$.
\end{prop}

\noindent{\bf Proof.} \ For each set $X$, let $F(X)= [[\Om^X],
\Om]^{\rm op}$. Then $F(X)$ is completely distributive. Let
$\eta_X: X\lra F(X)$ be given by $\eta_X(x)(\lam)=\lam(x)$ for all
$x\in X$ and $\lam\in\Om^X$.

It suffices to show that for any completely distributive
$\Om$-lattice $A$ and any function $f: X\lra A$, there is a unique
complete lattice morphism $g: F(X)\lra A$ such that $f=
g\circ\eta_X$.

{\bf Existence}: Since $f: X\lra A$, we have an $\Om$-functor
$f^\la:[A^{\rm op}, \Om]\lra[\Om^X]$, and then an $\Om$-functor
$(f^\la)^\la: [[\Om^X], \Om]^{\rm op}\lra[[A^{\rm op}, \Om],
\Om]^{\rm op}$.

Since $[A^{\rm op}, \Om]$ is completely distributive, $[A^{\rm
op}, \Om]^{\rm op}$ is also completely distributive by Theorem
1.1. Therefore, the $\Om$-functor $\inf_{[A^{\rm op},
\Om]}:[[A^{\rm op}, \Om], \Om]^{\rm op}\lra[A^{\rm op}, \Om]$ has
a right adjoint, hence it is complete lattice morphism. Let $g$ be
the composition of the following functors:
$$[[\Om^X], \Om]^{\rm
op}\stackrel{(f^\la)^\la}{\lra}[[A^{\rm op}, \Om], \Om]^{\rm
op}\stackrel{\inf_{[A^{\rm op}, \Om]}}{\lra}[A^{\rm op}, \Om]
\stackrel{\sup_A}{\lra} A.$$

Then $g$ is a complete lattice morphism. It remains to show that
$f=g\circ\eta_X$. Indeed, for each $x\in X$ and $\lam\in[A^{\rm
op}, \Om]$,
$$(f^\la)^\la)(\eta_X(x)(\lam)=\eta_X(x)(f^\la(\lam))=\lam(f(x)).$$
Thus, \begin{eqnarray*}\inf(f^\la)^\la(\eta_X(x)) &=&
\bw_{\lam\in[A^{\rm op},
\Om]}\Big((f^\la)^\la(\eta_X(x))(\lam)\ra\lam\Big)\\
&=&\bw_{\lam\in[A^{\rm op}, \Om]}\Big(\lam(f(x))\ra\lam\Big)\\
&=& \y(f(x)).
\end{eqnarray*}
Consequently, $g\circ\eta_X(x) = \sup\y(f(x))=f(x)$.

{\bf Uniqueness}: At first, we show that for all $G\in[[\Om^X],
\Om]$,
$$G=\bv_{\lam\in\Om^X}\Big[G(\lam)*\Big(\bw_{x\in
X}(\lam(x)\ra\eta_X(x))\Big)\Big].$$ Indeed, for all
$\mu\in\Om^X$, $$\Big(\bw_{x\in X}(\lam(x)\ra\eta_X(x))\Big)(\mu)
= \bw_{x\in X}(\lam(x)\ra\mu(x))=[\Om^X](\lam, \mu).$$ Thus,
$$\Big[\bv_{\lam\in\Om^X}\Big(G(\lam)*\Big(\bw_{x\in
X}(\lam(x)\ra\eta_X(x))\Big)\Big)\Big](\mu)=
\bv_{\lam\in\Om^X}G(\lam)*[\Om^X](\lam, \mu)=G(\mu).$$

Therefore, in $[[\Om^X], \Om]^{\rm op}$,
$$G=\bw_{\lam\in\Om^X}\Big[G(\lam)\rightarrowtail\Big(\bv_{x\in
X}(\lam(x)\otimes \eta_X(x))\Big)\Big].$$

Suppose $g:F(X)\lra A$ is a complete lattice morphism with
$f=g\circ \eta_X$. Then
$$g(G)=\bw_{\lam\in\Om^X}\Big[G(\lam)\rightarrowtail\Big(\bv_{x\in
X}(\lam(x)\otimes f(x))\Big)\Big].$$ Consequently, $g$ is unique.\
\ $\Box$

\end{document}